\magnification=1000
\hsize=11.7cm
\vsize=18.9cm
\lineskip2pt \lineskiplimit2pt
\nopagenumbers

\hoffset=-1truein
\voffset=-1truein

\advance\voffset by 4truecm
\advance\hoffset by 4.5truecm

\newif\ifentete

\headline{\ifentete\ifodd	\count0 
      \rlap{\head}\hfill\tenrm\llap{\the\count0}\relax
    \else
        \tenrm\rlap{\the\count0}\hfill\llap{\head} \relax
    \fi\else
\global\entetetrue\fi}

\def\entete#1{\entetefalse\gdef\head{#1}} 
\entete{}

\input amssym.def
\input amssym.tex

\def\-{\hbox{-}}
\def\.{{\cdot}}
\def\O{{\cal O}}

\def\F{{\cal F}}

\def\L{{\cal L}}

\def\N{{\cal N}}
\def\P{{\cal P}}

\def\G{{\cal G}}
\def\T{{\cal T}}

\def\I{{\cal I}}

\def\H{{\cal H}}

\def\A{{\cal A}}
\def\B{{\cal B}}

\def\C{{\cal C}}

\def\ad{\frak a\frak c}

\def\Ab{\frak A\frak b}

\def\Fct{\frak F\frak c\frak t}

\def\id{\frak i\frak d}
\def\int{\frak i\frak n\frak t}

\def\qq{\quad{\rm and}\quad}

\def\mod{\frak m\frak o\frak d}

\def\too{\longrightarrow}

\def\Set{\frak S\frak e\frak t}

 3
 2
\font\large=cmr10  scaled \magstep 2
 2
 2
 1
 2
\font\cds=cmr7
\font\cdt=cmti7

\count0=1

\centerline{\large A criterion on trivial homotopy}
\smallskip
\centerline{\bf Lluis Puig }
\par
\noindent 
\centerline{\cds CNRS, Institut de Math\'ematiques de Jussieu, lluis.puig@imj-prg.fr}
\par
\noindent
\centerline{\cds 6 Av Bizet, 94340 Joinville-le-Pont, France}

\medskip
\noindent
{\bf Abstract:} {\cds   In {\cdt Homotopy
decomposition of classifying spaces via elementary Abelian subgroups\/}, Stephan Jackowski and James McClure show, for functors admitting a {\cdt Mackey complement\/}  over categories with  {\cdt direct products\/}, a general result on vanishing cohomology. We~develop a  framework leading to a general result on {\cdt trivial homotopy\/} which partially generalizes Jackowski and McClure's result in two different directions.\/}

\bigskip
\noindent
{\bf £1. Introduction }
\medskip
£1.1.  In [2,~Corollary~5.16], Stephan Jackowski and James McClure state a quite general criterion on categories holding a {\it direct product\/}, which guaran-tees the annulation of the {\it positive cohomology groups\/} for the {\it contravariant $\Ab\-$valued\/} functors admitting a {\it Mackey complement\/}. In [3,~Theorem~6.26], we apply their argument to obtain a vanishing result  on the context of a {\it Frobenius $P\-$category $\F$\/} [3,~2.8] for a finite $p\-$group $P\,,$ since an example of such categories is the {\it additive cover\/} ---  introduced in [2] and recalled in section~£4 below ---  of the {\it full\/}  
subcategory~$\widetilde\F^{^{\rm sc}}$  of the 
{\it exterior quotient\/}~$\widetilde\F$~[3,~£1.3], over the set of  
{\it $\F\-$selfcentralizing\/} subgroups of~$P$~[3,~4.8];
denote by $\ad (\widetilde\F^{^{\rm sc}})$ and $\ad (\widetilde\F)$ their respective {\it additive covers\/}.

\medskip
£1.2.  But, we have~noticed a general kind of  {\it contravariant $\mod_\O\-$valued  $\ad (\widetilde\F^{^{\rm sc}})\-$functors\/} which need not  admit a {\it Mackey complement\/} and nevertheless have a  {\it  homotopically trivial standard  complex\/} (cf. defined in~£2.3 below). On the other hand, $\ad(\widetilde\F)$ have a family 
of commutative {\it square diagrams\/} (cf.~£5.1 below) --- which are not 
{\it pull-backs\/} --- in such a way that  any  {\it contravariant $\mod_\O\-$valued
  $\ad(\widetilde\F)\-$functor\/} which, relative to this family, admits a {\it complement\/} (cf.~£5.3) supplies again a  {\it  homotopically trivial standard complex\/}.

\medskip
£1.3. Our proofs of these two facts were similar enough to suggest the existence of a common argument leading to a more general result on {\it trivial homotopy\/}. The
effort in finding a common argument leads us to a somewhere sophisticated framework,
including all the situations described above; this general framework is described 
in section 2. The advantage of this point of view is that, once the {\it homotopy map\/} is defined,  the main result on {\it trivial homotopy\/} in section 3 consist on a simple
checking. The first application in section 4 includes  a general discussion relating the existence of a {\it direct product\/} in $\ad (\widetilde\F^{^{\rm sc}})$ with the fact that any $\widetilde\F^{^{\rm sc}}\!\-$morphism is an {\it epimorphism\/} [3,~Corollary~4.9]; the main result of this section is the key step in the proof of [4,~Theorem~6.22], as we explain in
£4.13 and £4.14 below. Similarly, in~£5.11 below we explain how the main result in 
section 5  can be applied in the proof of [4,~Theorem~8.10].

\vfill
\eject

\bigskip
\noindent
{\bf £2. The general framework}

\medskip
£2.1. We denote by $\Ab$ the category of  Abelian groups, by $\mod_\O$ the category of finitely generated $\O\-$modules where $\O$ is a complete discrete va-luation ring with unequal characteristics, and by $\Set$ the category of sets; all the other categories we consider are assumed to be {\it small\/}. For any cate-gory
$\C$ and any  pair of $\C\-$objects $Q$ and $R\,,$ we denote by $\C (Q,R)$ the set of $\C\-$morphisms from $R$ to $Q$ and set $\C (Q,Q) = \C (Q)$ for short; let us say that $\C$ is {\it ordered\/} if it fulfills the following condition [3,~A5.1.1]
\smallskip
\noindent
£2.1.1\quad{\it Any pair of $\C\-$objects $Q$ and $R$ fulfilling $\C (Q,R)\not= \emptyset\not= \C (R,Q)$ are $\C\-$isomorphic and then all the  
$\C\-$morphisms $R\to Q$ are  $\C\-$isomorphims.\/}
\smallskip
\noindent
In this case, note that $\C (Q)$ is a group.

\medskip
£2.2. Let $\A$ be an ordered category; we call {\it $\A\-$category\/} any
ordered category $\B$ containing $\A$ as a subcategory having the same objects, and fulfilling the following two conditions 
\smallskip
\noindent
£2.2.1.\quad {\it Any $\B\-$morphism $\varphi\,\colon R\to Q$ is the composition of a $\B\-$isomorphism 
$\varphi_*\,\colon R\cong R_*$ and an $\A\-$morphism $\iota\,\colon R_*\to Q\,.$\/}
\smallskip
\noindent
£2.2.2.\quad {\it Any $\B\-$isomorphism $\tau\,\colon R\cong R'$ such that $\iota'\circ \tau$  is an $\A\-$morphism
for some  $\A\-$morphism $\iota'\,\colon R'\to Q'\,,$ is an $\A\-$isomorphism.\/}
\smallskip
\noindent
In some sense, these conditions generalize the {\it divisibility\/} condition [3,~2.3] in~$\F$ when $\A = \F_P\,,$ the Frobenius $P\-$category of $P\,.$

\medskip
£2.3. Let $\B$ be an $\A\-$category and $\frak a\,\colon \B\to \Ab$ a 
{\it contrvariant\/} functor. In order to describe our {\it differential complex\/},
we fix a subcategory $\G$ of $\B$ having the same objects and only having 
$\G\-$isomorphisms.  For any  $n\in \Bbb N\,,$ setting 
$\Delta_n = \{i\mid 0\le i\le n\}$ considered as a category with the morphisms
$i\bullet i'$ defined by the ordered pairs $i\le i'\,,$ and denoting by 
$\Fct(\Delta_n,\B)$ the set of functors from $\Delta_n$ to $\B$ --- called 
{\it $\B\-$chains\/}   [3,~A2.2 and~A2.8] --- we set
 $$\Bbb C^n (\B, \frak a) = \prod_{\frak q\in 
 \Fct (\Delta_n,\B)} \frak a\big(\frak q(0)\big)
\eqno £2.3.1\phantom{.}$$ 
and consider the usual {\it differential map\/} sending $a = (a_\frak q)_{\frak q\in \Fct (\Delta_n,\B)}$ to the $\Bbb C^{n+1} (\B, \frak a)\-$element $d_n(a) = \big(d_n(a)_\frak r\big)_{\frak r\in 
\Fct (\Delta_{n+1},\B)}$ defined by [3,~A3.11.2]
$$d_n(a)_\frak r = \big(\frak a (\frak r (0\bullet 1))\big)(a_{\frak r\circ\delta_0^n})  + \sum_{i=1}^{n+1} (-1)^i a_{\frak r\circ\delta_i^n}
\eqno £2.3.2.$$
In this situation, we say that  $a = (a_\frak q)_{\frak q\in \Fct (\Delta_n,\B)}$ is 
{\it $\G\-$stable\/} [3,~A3.17] if, for any  {\it natural $\G\-$isomorphism\/} 
$\chi\,\colon \frak q\cong \frak q'$ between two {\it $\B\-$chains\/} $\frak q$ and $\frak q'$ in~$\Fct(\Delta_n,\B\,)\,,$  we have 
$a_{\frak q} =  \big(\frak a (\chi_0)\big)(a_{\frak q'})$
Then,  our {\it standard complex\/} is the {\it differential subcomplex\/} formed by
the subgroup of {\it $\G\-$stable\/} elements $\Bbb C^n_{\G} (\B, \frak a)$  of 
$\Bbb C^n (\B, \frak a)\,,$ and our main purpose is to give a criterion  
guaranteeing  that this {\it standard complex\/} admits a {\it homotopy map\/}.
\eject

\medskip
£2.4. Recall that an {\it interior structure\/} $\I$  in $\B$  is a correspondence sending any $\B\-$object $Q$ to a subgroup $\I (Q)$ of  $\B (Q)$ in such a way that we have
$$\varphi\circ\I (R)\i \I(Q)\circ \varphi
\eqno £2.4.1\phantom{.}$$
 for any $\B\-$morphism $\varphi\,\colon R\to Q$  [3,~1.3]; similarly, let us call 
 {\it co-interior structure\/}   $\I^\circ$  of $\B$ any {\it interior structure\/}  of
 the opposite  category  $\B^\circ\,;$ note that, according to 2.1, the correspondence  sending any $\A\-$object $Q$ to the intersection 
 $\I^\circ (Q)\cap \A(Q)$ is a {\it co-interior structure\/} of $\A\,.$
Then, a {\it bi-interior structure\/} $(\I,\I^\circ)$ in~$\B$ is any pair formed by an  {\it interior\/} and a {\it co-interior structures\/} $\I$ and $\I^\circ$ in $\B$ 
 such that $\I(Q) $ and $\I^\circ (Q)$ centralize each other for any  
 $\B\-$object~$Q\,;$ we  denote by~$\widetilde\B$ the corresponding {\it bi-exterior quotient\/} namely, for any pair of  $\B\-$objects $Q$ and $R\,,$ we set
 $$\widetilde\B (Q,R) = \I(Q)\backslash \B (Q,R)/\I^\circ (R)
 \eqno £2.4.2,$$
 and denote by $\frak e_{\I,\I^\circ}\,\colon \B\to \widetilde\B$ the canonical functor.

 \medskip
 £2.5. For any functor $\frak s\,\colon \B\to \Set$ we denote by 
 $\frak s\!\rtimes \B$ the new category --- called {\it semidirect product of $\B$ by 
 $\frak s$\/} --- where an object  is a pair $(s,Q)$ formed by a $\B\-$object $Q$ and an element $s\in \frak s (Q) = \frak s_Q\,,$  where a morphism $\varphi\,\colon (t,R)\to (s,Q)$  is just a  $\B\-$morphism 
 $\varphi\,\colon R\to Q$ such that $\frak s_\varphi (t) = s\,,$ and where the composition is determined by the composition in $\B$  [3,~A2.7]; in order to avoid confusion, we  denote this morphism by the pair $(t,\varphi)\,;$~more-over, we denote by $\frak s_\G$ the restriction of $\frak s$ to $\G$ and by 
 $\frak p\,\colon \frak s\!\rtimes \B\to \B$ the {\it forgetful\/} functor mapping $(s,Q)$ on~$Q\,.$

\medskip
£2.6. Our criterion to get a {\it homotopy map\/} is based on the existence of the following {\it data\/}. For any triple $(\A,\B,\G)$ as above, let us call {\it homotopic system\/} $\H = (\I,\I^\circ,\frak s,\frak n,\nu)$ any 
quintuple   formed by 
\smallskip
\noindent
£2.6.1\quad {\it A {\it bi-interior structure\/}  $(\I,\I^\circ)$ of $\B$  such that 
we have $\I (Q)\i \A(Q)$ and $\I (Q)\.\I^\circ (Q)\i \G (Q)$ for any 
$\B\-$object $Q\,,$ and that the corresponding bi-exterior quotient $\widetilde\A$ has a final object $P\,.$\/}

\smallskip
\noindent
£2.6.2\quad {\it A  functor $\frak s\,\colon \B\to \Set$  mapping any 
$\B\-$object~$Q$ on a finite set $\frak s_Q\,.$\/}
\smallskip
\noindent
£2.6.3\quad {\it A functor $\frak n \,\colon \frak s\rtimes\B\to \widetilde\A\i \widetilde\B$ sending $\frak s_\G\rtimes \G$ to $\frak e_{\I,\I^\circ}(\G) 
= \widetilde\G\,.$\/}
\smallskip
\noindent
£2.6.4\quad {\it A natural map $\nu \,\colon \frak n \to 
\frak e_{\I,\I^\circ}\circ\frak p = \tilde\frak p\,.$\/}
\smallskip
\noindent
Note that, for any $\B\-$object  $Q\,,$  any $s\in \frak s_Q$ and any element $\xi\in \I (Q)\.\I^\circ (Q)\,,$  the  
$\widetilde\A\-$morphism 
$$\frak n (s,\xi) : \frak n\big(s,Q\big)\cong \frak n\big(\frak s_\xi(s),Q\big)
\eqno £2.6.5\phantom{.}$$
 is a   $\widetilde\G\-$isomorphism.

\medskip
£2.7. In this situation, for any $n\in \Bbb N\,,$ any {\it $\B\-$chain\/} $\frak q\,\colon \Delta_n\to \B$ can  be  lifted,  for any choice of  
$s\in \frak s_{\frak q(0)}\,,$ to a unique {\it $\frak s\!\rtimes\B\-$chain\/}  
$\hat\frak q_s\,\colon \Delta_n\to \frak s\!\rtimes\B$ fulfilling 
$\hat\frak q_s (0) = \big(s,\frak q (0)\big)\,.$ Moreover, 
for any {\it $\frak s\!\rtimes\B\!\-$chain\/} $\hat\frak q\,\colon  \Delta_n\to \frak s\!\rtimes\B$  we have both the {\it $\widetilde\A\-$chain\/} 
$\frak n\circ \hat\frak q$ and, considering $\frak n\circ \hat\frak q$ as a {\it $\widetilde\B\-$chain\/}, the {\it natural\/} map
$$\nu *\hat\frak q : \frak n\circ\hat\frak q\too \tilde\frak p \circ\hat\frak q
\eqno £2.7.1;$$
then, for any $\ell\in \Delta_n\,,$  as in [3,~Lemma~A4.2] we denote by   
$$^{n} (\nu *\hat\frak q) :  {\Delta}_{n+1}\,\too \widetilde\B
\eqno £2.7.2\phantom{.}$$
 the functor which  coincides with $\frak n\circ\hat\frak q$ over~$\Delta_\ell\,,$  maps  $i\in \Delta_{n+1} -\Delta_\ell$  on $(\tilde\frak p\circ \hat\frak q) (i -1)\,,$ maps $i\bullet i+1$ on  $(\tilde\frak p \circ\hat\frak q )(i-1\bullet i)$ if $i\le n\,,$ and maps 
 $\ell\bullet \ell\!+\!1$  on 
 $$(\nu *\hat\frak q)_\ell : (\frak n\circ\hat\frak q) (\ell)\too 
 (\tilde\frak p\circ\hat\frak q) (\ell)
  \eqno £2.7.3;$$  
 moreover, denote by 
 $$\frak h_{n+1}^{n}(\nu *\hat\frak q) : \Delta_{n+1}\too 
 \widetilde\A \i \widetilde\B
 \eqno £2.7.4\phantom{.}$$
  the {\it $\widetilde\A\-$chain\/} $(\frak n\circ\hat\frak q)^P$ which extends 
$\frak n\circ \hat\frak q$ to $\Delta_{n+1}$ sending $n +1$ to $P$ and 
$n\bullet n +\!1$  to the unique $\widetilde\A\-$morphism from
 $(\frak n\circ\hat\frak q)(n)$ to $P\,.$

\medskip
£2.8. Then,  for any 
 {\it natural $\frak s_\G\!\rtimes\G\-$isomorphism\/}  $\hat\chi\,\colon \hat\frak q\cong \hat\frak q'$ we claim that we have the {\it natural $\widetilde\G\-$isomorphism\/}
$$^{n}  (\nu \diamond\hat\chi) : ^{n}  (\nu *\hat\frak q)\cong ^{n} (\nu *\hat\frak q')
\eqno £2.8.1\phantom{.}$$
which  coincides with $\frak n *\hat\chi$ over~$\Delta_\ell$ if $\ell\le n\,,$ and either  maps  $i\in \Delta_{n+1} -\Delta_\ell$  on $(\tilde\frak p *\hat\chi )_{i -1} $ or maps $\ell = n+1$
 on the {\it identity $\widetilde\B\-$morphism\/} of $P\,;$ indeed, for any  $0\le i\le \ell\le n\,,$ we have (cf. condition~£2.6.3)
$$^{n} (\nu \diamond\hat\frak q)(i) = 
\frak n \big(\hat\frak q(i)\big)\buildrel \frak n (\hat\chi_i)\over\cong 
\frak n \big(\hat\frak q'(i)\big) = ^{n} (\nu \diamond\hat\frak q')(i)
\eqno £2.8.2;$$
moreover, for any $\ell  +1 \le  i\le n+1 \,,$ we also have
$$^{n} (\nu \diamond\hat\frak q) (i ) 
= \tilde\frak p\big(\hat\frak q(i\! -1\!)\big)
\buildrel \tilde\frak p (\hat\chi_{i-1})\over\cong 
\tilde\frak p\big(\hat\frak q'(i\! -1\!)\big) = 
^{n} (\nu \diamond\hat\frak q') (i ) 
\eqno £2.8.3;$$
finally, if $\ell\le n$ then we get the commutative diagram  (cf. condition~£2.6.4)
$$\matrix{^{n} (\nu \diamond\hat\frak q) (\ell +1 ) 
= \tilde\frak p\big(\hat\frak q(\ell)\big)
&\buildrel \tilde\frak p (\hat\chi_\ell)\over\cong &
\tilde\frak p\big(\hat\frak q'(\ell)\big) = 
^{n} (\nu \diamond\hat\frak q') (\ell +1 )\cr
{\scriptstyle \nu_{\hat\frak q(\ell)}}\big\uparrow\hskip-60pt&\phantom{\Big\uparrow}
&\hskip-60pt\big\uparrow{\scriptstyle \nu_{\hat\frak q'(\ell)}}\cr
^{n} (\nu \diamond\hat\frak q) (\ell  ) = \frak n\big(\hat\frak q(\ell)\big)\hskip-20pt 
&\buildrel \frak n (\hat\chi_\ell)\over\cong &
\hskip-20pt \frak n\big(\hat\frak q'(\ell)\big) = ^{n} (\nu \diamond\hat\frak q') (\ell  )\cr} 
\eqno £2.8.4.$$

\medskip
£2.9. On the other hand,  any {\it contravariant\/} functor   
$\tilde\frak a\,\colon \widetilde\B\to \Ab$ determines a new  {\it contravariant\/} functor  from $\B$ to $\Ab$  sending any $\B\-$object $Q$ to the product 
 $\prod_{s\in \frak s_Q} (\tilde\frak a\circ\frak n)(s,Q)$ and any $\B\-$morphism $\varphi\,\colon R\to Q$ to the group homomorphism 
 $$\prod_{s\in \frak s_Q} (\tilde\frak a\circ\frak n)(s,Q)\too \prod_{t\in \frak s_R} (\tilde\frak a\circ\frak n)(t,R)
 \eqno £2.9.1\phantom{.}$$
 \eject
\noindent
mapping the element $a = \sum_{s\in \frak s_Q} a_s\,,$ where $a_s$ belongs to $(\tilde\frak a\circ\frak n)(s,Q)\,,$ on 
$$\sum_{t\in \frak s_R} \big((\tilde\frak a\circ \frak n)(t,\varphi)\big) 
(a_{\frak s_\varphi (t)})\in \prod_{t\in \frak s_R} (\tilde\frak a\circ\frak n)(t,R)
 \eqno £2.9.2;$$
 in particular, $\I(Q)$ acts on $\prod_{s\in \frak s_Q} (\tilde\frak a\circ\frak n)(s,Q)$ and we get a {\it contravariant\/} subfunctor
 $\H (\tilde\frak a)\,\colon \B\to \Ab$ sending any $\B\-$object $Q$ to
 $$\big(\H (\tilde\frak a)\big)(Q) = \big(\prod_{s\in \frak s_Q} (\tilde\frak a\circ\frak n)(s,Q)\big)^{\I(Q)}
 \eqno £2.9.3;$$
  indeed, for any $\rho\in \I (R)$ we have  $\varphi \circ \rho = \chi\circ\varphi$ for a suitable $\chi\in \I (Q)$ and therefore,
if $a = \sum_{s\in \frak s_Q} a_s$ belongs to $\big(\H (\tilde\frak a)\big)(Q)$ then for any $s\in \frak s_Q$ we have
$$\big((\tilde\frak a\circ \frak n)(s,\chi)\big) (a_{\frak s_\chi (s)}) = a_s
\eqno £2.9.4,$$
so that it is easily checked that
$$\eqalign{\Big(\big(\H (\tilde\frak a)\big)(\rho)\Big) 
\Big(\sum_{t\in \frak s_R} &\big((\tilde\frak a\circ \frak n)(t,\varphi)\big) (a_{\frak s_\varphi (t)})\Big)\cr
&= \sum_{t\in \frak s_R} \big((\tilde\frak a\circ \frak n)(t,\varphi)\big) (a_{\frak s_\varphi (t)})\cr }
 \eqno £2.9.5;$$
thus, the image of $a$ through homomorphism~£2.9.1 belongs to $\big(\H (\tilde\frak a)\big)(R)\,;$ moreover, the above correspondence sending $\tilde\frak a$ to $\tilde\H (\frak a)$  is clearly functorial.

\medskip
£2.10. Then, setting $\frak a = \tilde\frak a\circ \frak e_{\I,\I^\circ}\,,$   we claim that there is a {\it natural map\/}
$$\Delta_\H (\tilde\frak a) : \frak a \too \H(\tilde\frak a)
\eqno £2.10.1\phantom{.}$$
sending any $\B\-$object $Q$ to the group homomorphism
$$\Delta_\H (\tilde\frak a)_Q : \frak a (Q) \too \big(\prod_{s\in \frak s_Q} 
(\tilde\frak a\circ\frak n)(s,Q)\big)^{\I(Q)}
\eqno £2.10.2\phantom{.}$$
mapping $a\in \frak a (Q)$ on 
$$\Delta_\H (\tilde\frak a)_Q  (a) = \sum_{s\in \frak s_Q} (\tilde\frak a * \nu)_{(s,Q)}(a)
\eqno £2.10.3.$$
Indeed, this makes sense since $\tilde\frak p (s,Q)$ coincides with $Q\,;$ moreover, for  any $\B\-$morphism 
$\varphi\,\colon R\to Q\,,$ we have (cf.~condition~£2.6.4)
$$\eqalign{\Big(\big(\H (\tilde\frak a)\big)&(\varphi)\circ 
\Delta_\H (\tilde\frak a)_Q\Big)(a)
 = \sum_{t\in \frak s_R}  \big((\tilde\frak a\circ \frak n)(t,\varphi)\big) 
\big((\tilde\frak a * \nu)_{(\frak s_\varphi (t),Q)}(a)\big)\cr
& = \sum_{t\in \frak s_R} (\tilde\frak a * \nu)_{(t,R)}
\Big(\big(\frak a (\varphi)\big) (a)\Big)
= \big(\Delta_\H (\tilde\frak a)_R\circ \frak a (\varphi)\big)(a)\cr}
 \eqno £2.10.4.$$
We are ready to state our criterion on the  {\it contravariant\/} functor; we say that a {\it contravariant\/} functor $\tilde\frak a\,\colon \widetilde\B\to \Ab$ is 
{\it $\H\-$split\/} if the {\it natural map\/}  $\Delta_\H (\tilde\frak a)$ admits a  {\it natural section \/} 
$$\theta : \H (\tilde\frak a)\too \frak a
\eqno £2.10.5.$$
\eject

\bigskip
\noindent
{\bf £3. The trivial homotopy}

\medskip
£3.1. Let $\B$ be an $\A\-$category, $\G$ a subcategory of $\B$ having the same objects and only having $\G\-$isomorphisms (cf.~£2.3), $\H = (\I,\I^\circ,\frak s,\frak n,\nu)$ a {\it homotopic system\/} for the triple $(\A,\B,\G)$ (cf.~£2.6) and  $\tilde\frak a\,\colon \widetilde\B\to \Ab$  a  {\it  contravariant $\H\-$split\/} functor (cf.~£2.10), set $\frak a = \tilde\frak a\circ \frak e_{\I,\I^\circ}$ and denote by 
$$\theta : \H (\tilde\frak a)\too \frak a
\eqno £3.1.1\phantom{.}$$
a {\it natural section\/} of $\Delta_\H (\tilde\frak a)\,\colon \frak a \to 
\H(\tilde\frak a)$ (cf.~£2.10.1). Note that the group homomorphism sending  any element $(a_{\tilde\frak q})_{\tilde\frak q\in \Fct (\Delta_{n},\tilde\B)}$ of 
$\Bbb C^{n}_{\tilde\G} (\widetilde\B,\tilde\frak a)$ (cf.~£2.3) to the element 
$( a_{\frak q})_{\frak q\in \Fct (\Delta_{n},\B)}$ of 
$\Bbb C^{n}_{\G} (\B,\frak a)$ defined by $ a_\frak q = 
a_{\tilde\frak q}$ where $\tilde\frak q = \frak e_{\I,\I^\circ}\circ\frak q$ determines a group isomorphism
 $$\Bbb C^{n}_{\tilde\G} (\widetilde\B,\tilde\frak a)\cong \Bbb C^{n}_{\G} 
 (\B,\frak a)
 \eqno £3.1.2\phantom{.}$$
 and we identify to each other both members of this isomorphism. We denote by
 $\Bbb H^n_\G (\B,\frak a)$ the {\it $n\-$cohomology group\/} of our 
 {\it standard complex\/} (cf.~£2.3).

\medskip
£3.2. We are ready to define the {\it homotopy map\/};  for any $n\in \Bbb N\,,$ we consider the group homomorphism
$$ h^n : \Bbb C^{n+1}_{\G} (\B,\frak a)\too \Bbb C^n_{\G} (\B,\frak a)
\eqno £3.2.1\phantom{.}$$
sending any element $a = (a_{\frak r})_{\frak r\in \Fct (\Delta_{n+1},\B)}$ of 
$\Bbb C^{n+1}_{\G} (\B,\frak a)$ to the element $ h^n(a) =
\big( h^n( a)_{\frak q}\big)_{ \frak q\in \Fct (\Delta_n,\B)}$  of 
$\Bbb C^n_\G (\B,\frak a)$ defined by 
 $$ h^n( a)_{ \frak q} =  \sum_{\ell = 0}^{n+1} (-1)^\ell 
 \theta_{ \frak q(0)} \big(\sum_{s\in \frak s_{\frak q (0)} } 
a_{^n  ( \nu *\hat \frak q_s)}\big)
 \eqno £3.2.2.$$
 This makes sense since for any $\ell\in \Delta_{n+1}$  the element
 $b_\ell = \sum_{s\in \frak s_{\frak q (0)}} 
 a_{^n  ( \nu *\hat \frak q_s)}$
  belongs to $\big(\H (\tilde\frak a)\big)\big(\frak q (0)\big)\,,$ and moreover
  $h^n(a)$ is {\it $\G\-$stable\/}.
  
\medskip
£3.3.  Indeed, for any 
  $\xi_0\in \I\big(\frak q(0)\big)$ we can construct a  {\it natural $\B\-$auto-morphism\/}  $\xi\,\colon \frak q\cong \frak q$ which starts by $\xi_0\,\colon \frak q(0)\cong \frak q(0)$ and, for any 
$i$ in $\Delta_n - \{0\}\,,$ $\xi_i$ belongs to $\I\big(\frak q(i)\big)$  fulfilling
(cf.~2.4.1)
$$\frak q (i\! -\!1\bullet i)\circ \xi_{i-1} = \xi_i\circ \frak q (i\! -\!1\bullet i)
\eqno £3.3.1;$$
in particular, $\xi$ is a {\it natural  $\G\-$automorphism\/} (cf.~condition~£2.6.1) and therefore for any $s\in \frak s_{\frak q(0)}$ we get a  {\it natural $\frak s\!\rtimes \G\-$isomorphism\/} $\hat\xi_s : \hat\frak q_s\cong \hat\frak q_{\frak s_{\xi_0}(s)}$ (cf.~£2.7); then, from~£2.8 above,  for any $s\in \frak s_{\frak q(0)}$ and any 
$\ell\in \Delta_{n+1}$  we still get  a {\it natural $\widetilde\G\-$isomorphism\/}
$$^{n}  (\nu \diamond \hat\xi_s) : \frak h_\ell^{n}  (\nu *\hat\frak q_s)\cong \frak h_\ell^{n} (\nu *\hat\frak q_{\frak s_{\xi_0}(s)})
\eqno £3.3.2\phantom{.}$$
\eject
\noindent
and therefore the {\it $\G\-$stability\/} of $a = 
(a_{\frak r})_{\frak r\in \Fct (\Delta_{n+1},\B)}$  forces
$$ a_{\frak h_\ell^n  ( \nu *\hat \frak q_s)} =  \big((\tilde\frak a\circ \frak n) (\frak s_{\xi_0}(s),\xi_0)\big)( a_{\frak h_\ell^n  ( \nu *\hat \frak q_{\frak s_{\xi_0}(s)})})
\eqno £3.3.3;$$ 
at this point,  we have 
$$\xi_0 (b_\ell) = \sum_{s\in \frak s_{\frak q (0)}} \big((\tilde\frak a\circ \frak n) 
(s,\xi_0)\big) ( a_{\frak h_\ell^n  ( \nu *\hat \frak q_s)}) = 
\sum_{s\in \frak s_{\frak q (0)}} a_{\frak h_\ell^n  ( \nu *\hat \frak q_s)} = b_\ell
 \eqno £3.3.4.$$

 \medskip
£3.4. Moreover, for  any  {\it natural $\G\-$isomorphism\/}  
$\chi\,\colon \frak q\cong \frak q'$ and any $s\in \frak s_{\frak q(0)}$ 
we have a  {\it natural $\frak s_\G\rtimes\G\-$isomorphism\/} $\hat\chi_s : \hat\frak q_s\cong \hat\frak q'_{\frak s_{\chi_0}(s)}$ (cf.~£2.7);
 then, from~£2.8 above, we still get the $\frak s_\G\!\rtimes \G\-$isomorphism 
 $$\frak h_\ell^{n}  (\nu \diamond\hat\chi_s) : \frak h_\ell^{n} (\nu *\hat\frak q_s) \cong  \frak h_\ell^{n} (\nu *\hat\frak q'_{\frak s_{\chi_0}(s)})
\eqno £3.4.1\phantom{.}$$
and therefore the {\it $\widetilde\G\-$stability\/} of $a$ forces
 $$ a_{\frak h_\ell^{n}  (\nu *\hat\frak q_s)} = 
 \big((\tilde\frak a \circ\frak n) (\frak s_{\chi_0}(s),\chi_0)\big) 
 (a_{ \frak h_\ell^{n} (\nu *\hat\frak q'_{\frak s_{\chi_0}(s)})})
 \eqno £3.4.2.$$ 
 At this point,    setting $b'_\ell = \sum_{s'\in \frak s_{\frak q' (0)}}  
 a_{\frak h_\ell^n  ( \nu *\hat \frak q'_{s'})}$
for any $\ell\in \Delta_{n+1}\,,$ we have 
 $$ b'_\ell = \sum_{s\in \frak s_{\frak q (0)}}  \big((\tilde\frak a \circ\frak n) (\frak s_{\chi_0}(s),\chi_0)\big)^{-1}
( a_{\frak h_\ell^n  ( \nu *\hat \frak q_{s})}) = \big(\H(\tilde\frak a )\big) (\chi_0)^{-1}(b_\ell)
\eqno £3.4.3.$$
Consequently, we obtain
$$\eqalign{h^n( a)_{ \frak q'} 
&=  \sum_{\ell = 0}^{n+1} (-1)^\ell \theta_{ \frak q'(0)} (b'_\ell)
= \sum_{\ell = 0}^{n+1} (-1)^\ell \Big(\theta_{ \frak q'(0)} 
\circ\big(\H(\tilde\frak a )\big) (\chi_0)^{-1}\Big)(b_\ell)\cr
&= \sum_{\ell = 0}^{n+1} (-1)^\ell \big(\frak a (\chi_0)^{-1}\circ
\theta_{ \frak q(0)} \big)(b_\ell)
= \frak a (\chi_0)^{-1}\big(h^n(a)_\frak q\big)\cr}
\eqno £3.4.4.$$

\bigskip
\noindent
{\bf Theorem~£3.5.} {\it With the notation and the hypothesis above, for any $n\in \Bbb N$ we have
$$d_n\circ h^n  +  h^{n+1}\circ d_{n+1}  = 
{\rm id}_{\Bbb C^{n+1}_{\G} (\B,\frak a)}
\eqno £3.5.1.$$
In particular,  $\Bbb H^n_{\G} (\B,\frak a) = \{0\}$ for any $n\ge 1\,.$\/}

\medskip
\noindent
{\bf Proof:} For any $n\in \Bbb N\,,$ any 
$ a = ( a_{\frak r})_{\frak r\in \Fct (\Delta_{n+1},\B)}$ 
in $\Bbb C_{\G}^{n+1} (\B,\frak a)$ and any $\frak r\in \Fct (\Delta_{n+1},\B\,)\,,$ it suffices to prove that we have
$$ d_n\big( h^n ( a)\big)_{ \frak r} +  h^{n+1}\big( d_{n+1}( a)\big)_{\frak r}  =  a_{\frak r}
\eqno £3.5.2.$$
By the definition of $d_n$ (cf.~£2.3.2), for any {\it $\B\-$chain\/} $\frak r\,\colon\Delta_{n+1}\to \B$ we have
$$ d_n\big( h^n ( a)\big)_{\frak r} = \big( \frak a ( \frak r(0\bullet 1))\big)
\big( h^n ( a)_{ \frak r\circ\delta_0^n}\big) +
\sum_{i=1}^{n+1}(-1)^i h^{n} ( a)_{\frak r\circ\delta_i^n}
\eqno £3.5.3;$$
 \eject
\noindent
similarly, according to the very definition of $ h^{n+1}\,,$ we have
 $$ h^{n+1}\big( d_{n+1}( a)\big)_{\frak r}  =   \sum_{\ell = 0}^{n+2} (-1)^\ell   \theta_{ \frak r(0)} 
 \big(\sum_{s\in \frak s_{\frak r (0)}} d_{n+1} ( a)_{\frak h_\ell^{n+1}  (\nu *\hat \frak r_s)}\big)
 \eqno £3.5.4.$$

 \smallskip
But, for any $s\in \frak s_{\frak r(0)}$ and any $\ell\in \Delta_{n+2}$ we know that 
 $$\eqalign{ d_{n+1}( a)_{\frak h_\ell^{n+1}  (\nu *\hat \frak r_s)}\
 &= \Big( \tilde\frak a \big(\frak h_\ell^{n+1}  (\nu *\hat \frak r_s) (0\bullet 1)\big)\Big)
 ( a_{\frak h_\ell^{n+1}  (\nu *\hat \frak r_s) \circ\delta_0^{n+1}})\cr
 &\hskip50pt + \sum_{i=1}^{n+2} (-1)^i  a_{\frak h_\ell^{n+1}  (\nu *\hat \frak r_s)\circ\delta_i^{n+1}}\cr}
 \eqno £3.5.5;$$
 thus,  since $\frak h_0^{n+1}  (\nu *\hat \frak r_s) (0\bullet 1) = 
 \nu_{\hat \frak r_s (0)} = \nu_{(s,\frak r (0))}\,,$ from [3,~Lemma~A4.2]  it follows that  for $\ell = 0$ we get
  $$\eqalign{& d_{n+1}( a)_{\frak h_0^{n+1}  (\nu *\hat \frak r_s)} = 
  \big( \tilde\frak a(\nu_{(s,\frak r (0))})\big)  ( a_{ \frak r}) + 
  \sum_{i=1}^{n+2} (-1)^i  a_{\frak h_0^{n+1}(\nu *\hat\frak r_s)\circ
  \delta_i^{n+1}}\cr
 &= \big( \tilde\frak a (\nu_{(s,\frak r (0))})\big)( a_{\frak r}) -  a_{\frak h_0^{n+1}  (\nu *\hat\frak r_s)\circ\delta_1^{n+1}} - \sum_{i=1}^{n+1} (-1)^i   
 a_{\frak h_0^{n}  ( (\nu *\hat\frak r_s) *\delta_i^{n})}\cr}
 \eqno £3.5.6.$$
Moreover,  for any $1\le \ell\le n+2$ we have $\frak h_\ell^{n+1}  (\nu *\hat \frak r_s) (0\bullet 1) =
 ( \frak n\circ \hat\frak r_s)(0\bullet 1)\,;$ but,~according to~£2.6 above, 
 we denote the  $\frak s\!\rtimes \B\-$morphism 
$$\hat\frak r_s (0\bullet 1) : \big(s,\frak r (0)\big)\too 
\big(\frak s_{\frak r (0,1)}(s),\frak r (1)\big)
\eqno £3.5.7\phantom{.}$$
by the pair  $(s, \frak r (0\bullet 1))\,,$ and note that for any $i\in \Delta_{n+1}$ we have 
$$\eqalign{\frak h_{n+2}^{n+1}  (\nu *\hat\frak r_s)\circ\delta_i^{n+1} 
&=  (\frak n\circ\hat\frak r_s)^P\circ\delta_i^{n+1} = 
(\frak n\circ\hat\frak r_s\circ\delta_i^{n})^P\cr
& = \frak h_{n+ 1}^{n}  \big((\nu *\hat\frak r_s) *\delta_i^{n}\big)\cr}
\eqno £3.5.8;$$
consequently, from  [3,~Lemma~A4.2] again and from this equality, for any $1\le \ell\le n+2$ we obtain
 $$\eqalign{ d_{n+1}( a)_{\frak h_\ell^{n+1}  (\nu *\hat\frak r_s)}  
 &= \big( (\tilde\frak a \circ\frak n)(s,\frak r(0\bullet 1))\big) 
 ( a_{\frak h_{\ell-1}^{n}  ((\nu *\hat\frak r_s) *\delta_0^{n})})\cr
 &\hskip50pt + \sum_{i=1}^{n+2} (-1)^i  a_{\frak h_\ell^{n+1}  (\nu *\hat\frak r_s)\circ\delta_i^{n+1}}\cr}
 \eqno £3.5.9.$$

\smallskip
Similarly, it  follows  from [3,~Lemma~A4.2] that for any $1\le \ell\le n+2$ we still get
 $$\eqalign{\sum_{i=1}^{n+2} (-1)^i  & a_{\frak h_\ell^{n+1}  (\nu *\hat\frak r_s)\circ\delta_i^{n+1}} 
= \sum_{i=1}^{\ell -1} (-1)^i  a_{\frak h_{\ell-1}^n ((\nu *\hat\frak r_s) *\delta_i^n)}\cr
 &+ (-1)^\ell  a_{\frak h_\ell^{n+1}  (\nu *\hat\frak r_s)\circ\delta_\ell^{n+1}} 
 - (-1)^\ell  a_{\frak h_\ell^{n+1}  (\nu *\hat\frak r_s)\circ\delta_{\ell+1}^{n+1}}\cr
& - \sum_{i=\ell +1 }^{n+1} (-1)^i  a_{\frak h_\ell^{n} ((\nu *\hat\frak r_s) * \delta_i^n)}\cr}
\eqno £3.5.10,$$
 and from equality~£3.5.8 for $\ell = n+2$ we still obtain
$$\eqalign{\sum_{i=1}^{n+2} (-1)^i & a_{\frak h_{n+2}^{n+1}  (\nu *\hat\frak r_s)\circ\delta_i^{n+1}}\cr
& = \sum_{i=1}^{n+1} (-1)^i  a_{\frak h_{n+1}^n ((\nu *\hat\frak r_s) *\delta_i^n)}
 + (-1)^n  a_{\frak h_{n+2}^{n+1}  (\nu *\hat\frak r_s)\circ\delta_{n+2}^{n+1}} \cr}
\eqno £3.5.11.$$

\smallskip
At this point, for  any $s\in \frak s_{\frak r(0)}$ and any $1\le \ell\le n+2\,,$ 
let us set 
$$\eqalign{\alpha_{s,0} 
&= \big(\tilde\frak a (\nu_{(s,\frak r (0))})\big)( a_{\tilde\frak r})\cr
\alpha_{s,\ell} 
&= \Big((\tilde\frak a \circ\frak n)\big(s,\frak r(0\bullet 1)\big)\Big) 
 ( a_{\frak h_{\ell-1}^{n}  ((\nu *\hat\frak r_s) *\delta_0^{n})})\cr}
 \eqno £3.5.12;$$
 moreover, for any $0\le \ell\le n+1$ we still set [3,~Lemma~A4.2]
$$\eqalign{\beta_{s,\ell }
&= (-1)^\ell  a_{\frak h_\ell^{n+1}  (\nu *\hat\frak r_s)\circ\delta_{\ell +1}^{n+1}}
= (-1)^\ell  a_{\frak h_{\ell +1}^{n+1} (\nu *\hat\frak r_s)\circ\delta_{\ell+1}^{n+1}}\cr
0 &=\beta_{s,-1} = \beta_{s,n +2} = \gamma_{s,-1} = \gamma'_{s,n+2}\cr
\gamma_{s,\ell} 
&= \sum_{i=1}^{\ell } (-1)^i  a_{\frak h_{\ell}^n ((\nu *\hat\frak r_s) *\delta_i^n)}\qq
\gamma'_{s,\ell} 
= \sum_{i=\ell +1 }^{n+1} (-1)^i  a_{\frak h_\ell^{n} ((\nu *\hat\frak r_s) * \delta_i^n)}\cr}
\eqno £3.5.13.$$
With all this notation, from equalities~£3.5.6 and £3.5.9   for any~$s\in\frak s_{\frak r(0)}$  we get
$$\eqalign{\sum_{\ell = 0}^{n+2} (-1)^\ell  
 & d_{n+1}( a)_{\frak h_\ell^{n+1}  (\nu *\hat\frak r_s)}\cr
&= \sum_{\ell = 0}^{n+2} (-1)^\ell   (\alpha_{s,\ell } +\beta_{s,\ell -1} - \beta_{s,\ell } + \gamma_{s,\ell -1} - 
\gamma'_{s,\ell})\cr
&= \sum_{\ell = 0}^{n+2}(-1)^\ell  \alpha_{s,\ell}  - \sum_{\ell = 0}^{n +1}\sum_{i = 1}^{n +1}
(-1)^{\ell + i}  a_{\frak h_\ell^{n} ((\nu *\hat\frak r_s) * \delta_i^n)}\cr}
 \eqno £3.5.14.$$

 \smallskip
More precisely, for the first term of the bottom sum,  since $\tilde\frak a$ is {\it $\H\-$split\/},   we have (cf.~£3.5.12)
$$\theta_{ \frak r(0)}\big(\sum_{s\in \frak s_{\frak r (0)}}  \alpha_{s,0}\big)
= \theta_{\frak r(0)}\Big( \sum_{s\in \frak s_{\frak r (0)}}\big(\tilde\frak a (\nu_{(s,\frak r (0))})\big)( a_{\frak r})\Big) = a_{\tilde\frak r}
\eqno £3.5.15.$$
In conclusion, for any  $\ell\in \Delta_{n+1}$ and any $i\in \Delta_n$ setting 
$$\eqalign{\alpha_\ell &= \theta_{ \frak r(0)}\big(\sum_{s\in \frak s_{\frak r (0)}} \alpha_{s,\ell +1}\big)\qq
\varepsilon_{\ell,i} =  \theta_{ \frak r(0)}\big(\sum_{s\in \frak s_{\frak r (0)}} 
 a_{\frak h_\ell^{n} ((\nu *\hat\frak r_s) * \delta_{i +1}^n)}\big)\cr}
\eqno £3.5.16,$$
from £3.5.4,~£3.5.14 and £3.5.15 we get
$$ h^{n+1}\big( d_{n+1}( a)\big)_{\tilde\frak r} =  a_{\tilde\frak r} - \sum_{\ell = 0}^{n +1}(-1)^\ell 
\big(\alpha_\ell + \sum_{i = 1}^{n +1} (- 1)^i\varepsilon_{\ell,i-1}\big)
 \eqno £3.5.17.$$
 \eject

\smallskip
On the other hand, setting $\xi_\ell = \sum_{s\in \frak s_{\frak r( 0)}} 
a_{\frak h_\ell^n   (\nu *(\hat\frak r_s\circ\delta_0^n))}\,,$
  from the {\it naturality\/} of $\theta\,,$ for the first term in the right-hand member of equality~£3.5.3 we get
$$\eqalign{\big(\frak a (\tilde\frak r(0\bullet 1))\big) \big( h^n ( a)_{\frak r\circ\delta_0^n}\big)
& = \big(\frak a (\tilde\frak r(0\bullet 1))\big)\big(\sum_{\ell = 0}^{n+1} (-1)^\ell 
\theta_{\frak  r(1)} (\xi_\ell)\big)\cr
&= \sum_{\ell = 0}^{n+1} (-1)^\ell  \Big(\theta_{\frak r(0)}\circ \big(\H (\frak a)\big)\big( \frak r(0\bullet 1)\big)\Big) (\xi_\ell)\cr}
\eqno £3.5.18;$$
but, by the very definition of $\H(\tilde\frak a)$ we still get 
(cf.~£2.9.2 and~£3.5.12)
$$\eqalign{\Big(\big(\H (\tilde\frak a)\big)\big( \frak r(0\bullet 1)\big)\Big) (\xi_\ell)
& = \sum_{s\in \frak s_{\frak r(0)}} \big((\tilde\frak a\circ \frak n)(s, \frak r(0\bullet 1))\big) (a_{\frak h_\ell^n   (\nu *(\hat\frak r_s\circ\delta_0^n))})\cr
&= \sum_{s\in \frak s_{\frak r(0)}} \alpha_{s,\ell +1}\cr}
\eqno £3.5.19\,;$$
hence, we obtain (cf.~£3.5.16)
$$\big(\frak a (\frak r(0\bullet 1))\big) \big( h^n ( a)_{\frak r\circ\delta_0^n}\big) = \sum_{\ell = 0}^{n +1}(-1)^\ell \alpha_\ell
\eqno £3.5.20.$$

\smallskip
 Moreover, if $i\not= 0\,,$ for the $i\-$th  term in the right-hand member of equality~£3.5.3 and  for any 
 $s\in \frak s_{\frak r(0)}\,,$ it follows from~£2.7 that $(\widehat{\frak r\circ\delta_i^n})_s = \hat\frak r_s\circ \delta_i^n$
and therefore  we get
 $$ \eqalign{&h^n (a)_{\frak r\circ\delta_i^n} 
 = \sum_{\ell = 0}^{n+1} (-1)^\ell \theta_{\frak r(0)}  \big(\sum_{s\in \frak s_{\frak r (0)}} 
 a_{\frak h_\ell^n  (\nu *(\hat\frak r_s\circ\delta_i^n))}\big)\cr
& = \sum_{\ell = 0}^{n+1} (-1)^\ell \theta_{\frak r(0)}  \big(\sum_{s\in \frak s_{\frak r (0)}} 
a_{\frak h_\ell^n  ((\nu *\hat\frak r_s) *\delta_i^n)}\big)
 = \sum_{\ell = 0}^{n+1} (-1)^\ell \varepsilon_{\ell,i-1}\cr} 
  \eqno £3.5.21.$$
  Finally, equality~£3.5.2 follows from the suitable alternating sum of equalities~£3.5.17,~£3.5.20,  and~£3.5.21. We are done.

\bigskip
\noindent
{\bf £4. The direct product case}

\medskip
£4.1.  Let $\A$ be an ordered  category  and $\B$ an 
{\it  $\A\-$category\/} (cf.~£2.1); we assume that $\A$ and $\B$ have a
finite set of objects and that, for any $\B\-$object~$Q\,,$ the groups $\A (Q)$
and $\B (Q)$ are finite. Recall that the {\it additive cover\/}~$\ad (\B)$  [3,~A2.7.3] of~$\B$ --- introduced in~[2] --- is the category where the {\it objects\/} are the finite sequences  $\bigoplus_{i\in I} Q_i$ of $\B\-$objects $Q_i$ and where the 
{\it morphisms\/}
$$(f,\varphi ) : \bigoplus_{j\in J} R_j\too \bigoplus_{i\in I} Q_i
\eqno £4.1.1\phantom{.}$$
\eject
 \noindent
are the pairs formed by a map $f\,\colon J\to I$ and by a family $\varphi = \{\varphi_j\}_{j\in J}$ of $\B\-$morphisms
$$\varphi_j : R_j\too Q_{f(j)}
\eqno £4.1.2\phantom{.}$$
where $j$ runs over $J$, the {\it composition\/} being induced by the composition of maps and the composition in $\B\,.$ This category admits an obvious {\it direct sum\/} and we denote by 
$$\frak j : \B\too \ad (\B)
\eqno £4.1.3\phantom{.}$$
the canonical functor sending any $\B\-$object $Q$ to the
$\ad (\B)\-$object $\bigoplus_{\{\emptyset\}} Q$ that we still denote by $Q\,;$ 
note that any {\it contravariant\/} functor $\frak a\,\colon \B\too \Ab$ determines an {\it additive contravariant\/} functor
$$\frak a^{\ad} : \ad(\B)\too \Ab
\eqno £4.1.4.$$

\medskip
£4.2. Let us say that $\B$ is a {\it multiplicative $\A\-$category\/} whenever the 
{\it additive cover\/} $\ad (\B)$ has {\it direct products\/}. As we mention in the Introduction, one case we are interested in is the following. Let $P$ be a finite $p\-$group and $\F$ a Frobenius $P\-$category [3,~2.8], and denote 
by~$\widetilde\F$ the {\it exterior quotient\/} of~$\F$ [3,~£1.3] ---
that is to say, the quotient of $\F$ by the {\it inner\/} automorphisms of the 
$\F\-$objects --- and by $\widetilde\F^{^{\rm sc}}$ the {\it full\/} subcategory of 
$\widetilde\F$  over the set of  {\it $\F\-$selfcentralizing\/} subgroups of $P$ [3,~4.8]. We show in   [3,~Proposition~6.14] that $\ad (\widetilde\F^{^{\rm sc}})$ admits  a {\it direct product\/};
 consequently, denoting by $\F_P$ the Frobenius category of the group $P$ [3,~1.8] and by $\widetilde\F_P^{^{\rm sc}}$ the {\it full\/} subcategory of $\widetilde\F_P$  over the same set of  {\it $\F\-$selfcentralizing\/} subgroups,  $P$ is a final  $\widetilde\F_P^{^{\rm sc}}\-$object  and   $\widetilde\F^{^{\rm sc}}$ is clearly a {\it $\widetilde\F_P^{^{\rm sc}}\-$category\/} which  is {\it ordered\/} and {\it multiplicative\/}; thus,  Proposition~£4.12 in this section can be applied to   
 $\tilde\F^{^{\rm sc}}\,.$
  
\bigskip
\noindent
{\bf Lemma~£4.3.} {\it If the $\A\-$category $\B$ is multiplicative then any morphism is an epimorphism.\/} 

\medskip
\noindent
{\bf Proof:} For any $\B\-$object $Q$ and any $\B\-$automorphism 
$\sigma\in \B (Q)$ we have a unique $\ad (\B)\-$morphism 
$\delta_\sigma\,\colon Q\to Q\times Q$ which composed with the structural 
$\ad (\B)\-$morphisms of $Q\times Q$ yields ${\rm id}_Q$ and $\sigma\,;$ since $\B$ is an {\it ordered\/} category, up to suitable identifications, we have 
$$Q\times Q = \big(\bigoplus_{\sigma\in \B (Q)} Q\big)\oplus W
\eqno £4.3.1\phantom{.}$$
where $W$ is an $\ad (\B)\-$object such that there is {\it no\/} $\B\-$morphism from $Q$ to $W\,,$
and where the structural $\ad (\B)\-$morphisms
$$Q\longleftarrow \big(\bigoplus_{\sigma\in \B (Q)} Q\big)\oplus W \too Q
\eqno £4.3.2\phantom{.}$$
respectively induce ${\rm id}_Q\,,$ and $\sigma$ over the term indexed by 
$\sigma\in \B (Q)\,.$
\eject

\smallskip
 Let $\psi\,\colon T\to R\,,$  $\varphi\,\colon R\to Q$ and 
 $\varphi'\,\colon R\to Q$ be $\B\-$morphisms such that $\varphi'\circ\psi = \varphi\circ\psi \,;$ then, there is a unique $\ad (\B)\-$morphism 
 $\theta\,\colon R\to Q\times Q$ which composed with the structural 
 $\ad (\B)\-$morphisms of $Q\times Q$ yields $\varphi\,\colon R\to Q$
 and $\varphi'\,\colon R\to Q\,;$ moreover, the composition $\theta\circ\psi$ has to coincide with the unique
$\ad (\B)\-$morphism  
$$T\too Q\times Q = \big(\bigoplus_{\sigma\in \B (Q)} Q\big)\oplus W
\eqno £4.3.3\phantom{.}$$
which composed with the structural $\ad (\B)\-$morphisms of $Q\times Q$ yields
$\varphi\circ\psi $ and~$\varphi'\circ\psi\,;$ but, since $\varphi\circ\psi =\varphi'\circ\psi\,,$ this $\ad (\B)\-$morphism  necessarily~coin-cides with the $\ad (\B)\-$morphism formed by the map sending to ${\rm id}_Q\,,$ in the set of indices 
associated with $\big(\bigoplus_{\sigma\in \B (Q)} Q\big)\oplus W\,,$ the unique element $\emptyset$ in the set associated with $T\,,$ and by the $\B\-$morphism 
$\varphi\circ\psi \,\colon T\to Q$ (cf.~£4.1). Consequently,  the $\ad (\B)\-$morphism $\theta$ is also 
formed by  the map sending $\emptyset$ to~${\rm id}_Q$ and by a $\B\-$morphism $\varphi''\,\colon R\to Q$ which composed with the
identity map over $Q$ yields $\varphi$ and $\varphi'$ which forces
$\varphi = \varphi'' = \varphi'\,.$

 \bigskip
 \noindent
 {\bf Lemma~£4.4.} {\it If the $\A\-$category $\B$ is multiplicative then the
 {\it additive cover\/} $\ad (\B)$ admits pull-backs.\/}
 
\medskip
\noindent
{\bf Proof:} It is easily checked that it suffices to prove that any pair of 
$\B\-$mor-phisms $\alpha\,\colon Q\to T$ and $\beta\,\colon R\to T$ admits a
{\it pull-back\/}; but, if $U = \bigoplus_{\ell\in L} U_\ell$ is an $\ad (\B)\-$object
and 
$$\gamma = \{\gamma_\ell\}_{\ell\in L} : \bigoplus_{\ell\in L} U_\ell = U\to Q\qq \delta = \{\delta_\ell\}_{\ell\in L} : \bigoplus_{\ell\in L} U_\ell = U\to R
\eqno £4.4.1\phantom{.}$$ 
are a pair of  $\B\-$morphisms fulfilling $\alpha\circ\gamma  = \beta\circ\delta\,,$
we certainly have a unique $\ad (\B)\-$morphism 
$$\varepsilon : U\too Q\times R = \bigoplus_{i\in I} S_i
\eqno £4.4.2\phantom{.}$$
which composed with the structural $\ad (\B)\-$morphisms
$$\{\chi_i\}_{i\in I} : \bigoplus_{i\in I} S_i = Q\times R \to Q\qq \{\rho_i\}_{i\in I} : \bigoplus_{i\in I} S_i = Q\times R \to R
\eqno £4.4.3\phantom{.}$$
respectively yields $\gamma$ and $\delta\,.$

\smallskip
 That is to say, we have a map $e\,\colon L\to I$ and, for any $\ell\in L\,,$
a $\B\-$morphism $\varepsilon_\ell\,\colon U_\ell\to S_{e(\ell)}$ fulfllling
$$\chi_{e (\ell)} \circ \varepsilon_\ell = \gamma_\ell \qq
\rho_{e (\ell)} \circ \varepsilon_\ell  = \delta_\ell
\eqno £4.4.4\phantom{.}$$
and therefore still fulfilling
$$\alpha\circ \chi_{e (\ell)} \circ \varepsilon_\ell = \alpha\circ\gamma_\ell 
= \beta\circ \delta_\ell = \beta\circ \rho_{e (\ell)} \circ \varepsilon_\ell 
\eqno £4.4.5;$$
hence, it follows from Lemma~£4.3 that $\alpha\circ \chi_{e (\ell)} = \beta\circ \rho_{e (\ell)}$ for any $\ell\in L\,.$ Thus, denoting by $I_{\alpha,\beta}\i I$
the subset of $i\in I$ such that $\alpha\circ \chi_i = \beta\circ \rho_i\,,$
and setting $Q {}_\alpha\!\times_\beta R = \bigoplus_{i\in I_{\alpha,\beta}} S_i\,,$
in $\ad (\B)$ we obtain the {\it pull-back}
$$\matrix{&&Q {}_\alpha\!\times_\beta R\cr
&\swarrow\hskip-15pt&&\hskip-15pt\searrow\cr
Q\hskip-10pt&&&&\hskip-10pt R\cr
&{\atop \alpha}\hskip-5pt\searrow\hskip-15pt&&\hskip-15pt\swarrow\hskip-5pt{\atop \beta}\cr
&&T\cr}
\eqno £4.4.6.$$

\medskip
£4.5. As a matter of fact,  for any triple of $\B\-$objects $Q\,,$ $R$ and~$T\,,$ 
 if~$\B$ is an  $\A\-$category  where  any morphism is an epimorphism,  
any  $\B\-$morphism\/ $\alpha\,\colon Q\to R$ induces an injective map from 
$\B (T,R)$ to~$\B (T,Q)$ and, as in~[3,~6.4], we can consider the elements 
 of~$\B (T,Q)$ which, even {\it partially,\/} cannot be extended {\it via\/} 
 $ \alpha\,;$ precisely, we set
$$\B (T,Q)_{\alpha} = \B (T,Q) - \bigcup_{\theta'} \B (T,Q')\circ \theta'
\eqno £4.5.1,$$
where $\theta'$ runs over the set of {\it $\B\-$nonisomorphisms\/}
$\theta'\,\colon Q\to  Q'$ from $Q$  fulfilling $\alpha'\circ \theta' = \alpha$ for some $\alpha'\in \B (R,Q')$ which is then unique, and we simply 
 say that $\theta'$ {\it divides\/} $\alpha$ setting $\alpha' = \alpha/\theta'\,.$ Note that a $\B\-$morphism $\beta\,\colon Q\to T$ belongs to $\B (T,Q)_{\alpha}$ if and only if $\alpha$ belongs to $\B(R,Q)_\beta$ and,   
  as in [3,~6.9], let us call {\it strict\/} such a triple $(\alpha,Q,\beta)$ in  $\B\,;$ moreover, we say that two   {\it strict triples\/} $(\alpha,Q,\beta)$ and 
$(\alpha',Q',\beta')$ in $\B$ are {\it equivalent\/} if there exists a 
$\B\-$isomorphism $\theta\,\colon Q\cong Q'$ fulfilling
 $$\alpha'\circ\theta = \alpha\qq \beta'\circ\theta = \beta
 \eqno £4.5.2;$$
 let us denote by  $\frak T_{R,T}$  a set of representatives for the set of equivalence classes of  strict triples to $R$ and $T$  in $\B\,;$ note that our finiteness
 hypothesis in~£4.1 forces $\frak T_{R,T}$ to be finite. Then the analogous of equality [3,~6.7.1] still holds in this general setting and, in some sense, it characterizes   the {\it multiplicativeness\/} of $\B$ and the {\it direct product\/} in $\ad (\B)\,.$

\bigskip
\noindent
{\bf Proposition~£4.6.} {\it With the notation above, $\B$ is multiplicative if and only if any $\B\-$morphism is an epimorphism and, for any triple of $\B\-$objects $Q\,,$ $R$ and~$T\,,$ and any $\alpha\in \B (R,Q)\,,$ we have
$$\B (T,Q) = \bigsqcup_{ \theta'} \B (T,Q')_{\alpha/\theta'} \circ \theta'
\eqno £4.6.1\phantom{.}$$
where $\theta'\,\colon Q\to Q'$ runs over a set of representatives for the set of 
$\B\-$isomor-phism classes of $\B\-$morphisms starting on $Q$ and dividing 
$\alpha\,.$  In this case, we have a unique  $\ad(\B)\-$isomorphism
$$\bigoplus_{(\alpha',Q',\beta')\in \frak T_{R,T}} Q' \cong R\times T
\eqno £4.6.2\phantom{.}$$
which composed with the structural $\ad(\B)\-$morphisms from $R\times T$ to $R$ and $T$ respectively yields $\alpha'$ and $\beta'\,.$\/}
\eject

\medskip
\noindent
{\bf Proof:} Assume first that $\B$ is {\it multiplicative\/}; then, it follows from Lemma~£4.3 that  any $\B\-$morphism is an epimorphism. On the other hand,  once we fix $\alpha$ in~$\B (R,Q)\,,$ we have a bijection between $\B (T,Q)$ and the set of $\ad(\B)\-$mor-phisms from $Q$ to the {\it direct product\/} $R\times T = \bigoplus_{i\in I}\, Q_i$
which composed with the structural $\ad(\B)\-$morphism $R\times T\to R$ coincide with~$\alpha\,;$ moreover, the $\ad(\B)\-$morphism corresponding to 
$\beta\in \B(T,Q)$ is given by some $i\in I$ and by a suitable $\B\-$morphism  
$\theta_i\,\colon Q\to Q_i$ dividing $\alpha$ and $\beta\,.$

\smallskip
In particular, if $\beta$ belongs to $\B(T,Q)_\alpha$ then $\theta_i$ has to be a 
$\B\-$isomorphism. In any case, we have $\alpha = \alpha_i\circ \theta_i$ and 
$\beta = \beta_i\circ \theta_i$ for  suitable $\alpha_i\in \B (R,Q_i)$ and  
$\beta_i\in \B (T,Q_i)\,,$ and it is quite clear that $\beta_i$ necessarily belongs to 
$\B(T,Q_i)_{\alpha_i}\,;$ now, equality~£4.6.1 follows easily. Moreover, this argument applied to any $(\alpha',Q',\beta')\in \frak T_{R,T}$ determines a unique $i\in I$ together with a $\B\-$isomorphism $Q'\cong Q_i\,,$ which proves isomorphism~£4.6.2.

\smallskip
Conversely, assume that any $\B\-$morphism is an epimorphism and that equalities~£4.6.1 hold; it is easily checked that in order to exhibit  
{\it  direct products\/} in~$\ad (\B)\,,$ it suffices   to exhibit the 
{\it direct product\/} of any pair of $\B\-$objects $R$ and $T\,;$ thus,  choosing  
 a set of representatives $\frak T_{R,T}$ for the set of equivalence classes of  strict triples to $R$ and $T$ in $\B\,,$ we simply claim that 
 $$R\times T = \bigoplus_{(\alpha',Q',\beta')\in \frak T_{R,T}} Q'
 \eqno £4.6.3.$$
 Indeed, for any $\B\-$object $Q$ and any pair of $\B\-$morphisms $\alpha \,\colon Q\to R$ and $\beta\,\colon Q\to T$
 we have the corresponding equality~£4.6.1 determined by $Q\,,$ $R\,,$ $T$ and 
 $\alpha\,;$ then, the $\B\-$morphism $\beta$ determines a unique  term $\theta''\,\colon Q\to Q''$ in the disjoint union,
  together with a unique $\beta''\in \B (T,Q'')_{\alpha/\theta''}\,;$ that is to say, we get an $\ad (\B)\-$morphism
$$Q\too \bigoplus_{(\alpha',Q',\beta')\in \frak T_{R,T}} Q'
\eqno £4.6.4\phantom{.}$$
which composed with the structural homomorphisms yields $\alpha$ and $\beta\,,$ and the uniqueness is clear. We are done.

\bigskip
\noindent
{\bf Corollary~£4.7.} {\it With the notation above, if $\B$ is multiplicative
and $\I$ is an interior structure of $\B$ such that $\I (Q)$ is contained in $\A (Q)$  for any $\A\-$object~$Q\,,$ then the exterior quotient $\widetilde\B$ is multiplicative too, and the strict $\widetilde\B\-$triples are the image
of the strict $\B\-$triples.\/}

\medskip
\noindent
{\bf Proof:} If $\tilde\psi\,\colon T\to R\,,$  $\tilde\varphi\,\colon R\to Q$ and 
 $\tilde\varphi'\,\colon R\to Q$ are $\widetilde\B\-$morphisms fulfilling 
 $\tilde\varphi'\circ \tilde\psi = \tilde\varphi\circ \tilde\psi $ then, choosing respective representatives $\psi\in \B (R,T)$ and $\varphi\,,\,\varphi'\in \B (Q,R)\,,$
we have  $\varphi'\circ \psi = \theta\circ\varphi\circ \psi $ for some 
$\theta\in \I (Q)$ and therefore, according to Lemma~£4.3, we still have
$\varphi'= \theta\circ\varphi \,,$ so that $\tilde\varphi' = \tilde\varphi\,;$
that is to say, any $\widetilde\B\-$morphism is an epimorphism.

\smallskip
Moreover,  for any triple of $\widetilde\B\-$objects $Q\,,$ $R$ and~$T\,,$ and
any  $\widetilde\B\-$morphism\/ $\tilde\alpha\,\colon Q\to R\,,$ we claim that
$\widetilde\B (Q,T)_{\tilde\alpha}$ is the image in $\widetilde\B (Q,T)$ of 
$\B (Q,T)_{\alpha}\,;$ indeed, by the very definition~£4.5.1, it is clear that this image contains $\widetilde\B (Q,T)_{\tilde\alpha}\,;$  but, if $\tilde\beta$ belongs to $\widetilde\B (Q,T) - \widetilde\B (Q,T)_{\tilde\alpha}$ then there is a  
{\it $\widetilde\B\-$noniso-morphism\/} $\tilde\theta'\,\colon Q\to  Q'$ from 
$Q$  fulfilling $\tilde\beta = \tilde\beta'\circ \tilde\theta'$ and 
$\tilde\alpha =\tilde\alpha'\circ \tilde\theta'$ for some $\tilde\beta'\in 
\widetilde\B (T,Q')$ and  $\tilde\alpha'\in \widetilde\B (R,Q')\,;$ thus, for any
representatives $\beta$ of $\tilde\beta$ and $\alpha$ of $\tilde\alpha\,,$
 choosing respective representatives $\theta'\in \B (Q',Q)\,,$ $\beta'\in \B (T,Q')$  and $\alpha'\in \B (R,Q')\,,$ for suitable $\tau'\in \I (T)$ and $\rho'\in \I (R)$ we get 
$$\beta = \tau'\circ \beta'\circ \theta'\qq \alpha =\rho'\circ \alpha'\circ \theta'
\eqno £4.7.1,$$
so that $\beta$ belongs to $\B (Q,T) - \B (Q,T)_{\alpha}\,.$

\smallskip
Consequently, it follows from equality~£4.6.1 applied to $\B$ that
$$\widetilde\B (T,Q) = \bigcup_{ \tilde\theta'} \widetilde\B (T,Q')_{\tilde\alpha/\tilde\theta'} \circ \tilde\theta'
\eqno £4.7.2\phantom{.}$$
where $\tilde\theta'\,\colon Q\to Q'$ runs over a set of representatives for the set of $\widetilde\B\-$isomor-phism classes of $\widetilde\B\-$morphisms starting on $Q$ and dividing $\tilde\alpha\,;$ moreover, since $\I (T)$ acts by composition on the left stabilizing the term $\widetilde\B (T,Q')_{\tilde\alpha/\tilde\theta'} \circ \tilde\theta'$ for any~$\tilde\theta'$ this union is disjoint too. Now the corollary follows from Proposition~£4.6.

\medskip
£4.8. From now on, let  $\G$ be a subcategory   of $\B$  having  the same objects and only having   $\G\-$isomorphisms, and $\I$  an {\it interior structure\/} of $\B$ such that $\I (Q)$ is contained in $\A (Q)$ and $\G (Q)
$ for any $\A\-$object $Q\,,$ and that the {\it exterior quotient\/} $\widetilde\A$ has a final object $P\,.$ Moreover, we assume that
\smallskip
\noindent
£4.8.1\quad {\it The $\widetilde\A\-$category $\widetilde\B$ is multiplicative\/}. 
\smallskip
\noindent
Then, we have the functor
$$\widetilde{\frak m}_P : \widetilde\B\too \ad (\widetilde\B)
\eqno £4.8.2\phantom{.}$$
sending any $\widetilde\B\-$object $Q$ to the direct product $Q\times P$ and any 
$\widetilde\B\-$morphism $\tilde\varphi\,\colon R\to Q$ to the unique 
$\ad (\widetilde\B)\-$morphism
$$\widetilde\varphi\times \widetilde{\rm id}_P : R\times P\too Q\times P
\eqno £4.8.3\phantom{.}$$ 
determined by the pair of $\widetilde\B\-$morphisms
$$R\times P\too R\buildrel \tilde\varphi\over\too  Q\qq R\times P\too P\buildrel \widetilde{\rm id}_P\over\too P
\eqno £4.8.4;$$
note that we have a {\it natural map\/} $\tilde\omega_P \,\colon 
\widetilde{\frak m}_P \to \tilde\frak j$
sending any $\widetilde\B\-$object $Q$ to the structural $\ad(\widetilde\B)\-$morphism 
$Q\times P\to Q\,.$

\medskip
£4.9.   Explicitly, for any $\widetilde\B\-$object $Q\,,$ setting 
$\tilde I_Q = \frak T_{Q,P}$ we certainly may assume that
$$Q\times P = \bigoplus_{(\tilde\alpha',Q',\tilde\iota')\in \tilde I_Q} Q'
\eqno £4.9.1;$$
moreover, it follows from condition~£2.2.1 that, up to $\ad (\widetilde\B)\-$isomorphisms, we still may assume that $\tilde\iota' $ coincides with the unique
$\widetilde\A\-$morphism $\tilde\iota_{Q'}\to P$ and that the {\it strict triple\/} 
$(\widetilde{\rm id}_Q,Q,\tilde\iota_Q)$ belongs to~$\tilde I_Q\,.$ On the other hand,  for any $\widetilde\B\-$morphism $\tilde\varphi\,\colon R\to Q\,,$ the 
$\ad(\widetilde\B)\-$morphism $\tilde\varphi\times \widetilde{\rm id}_P$ determines both a map $\tilde I_{\tilde\varphi}\,\colon I_R\to I_Q$ and,  for any triple $(\tilde\beta',R',\tilde\iota_{R'})$ belonging to $\tilde I_R\,,$ setting 
$I_{\tilde\varphi}(\tilde\beta',R',\tilde\iota_{R'}) = 
(\tilde\alpha',Q',\tilde\iota_{Q'})$ a $\widetilde\B\-$morphism
$\tilde\varphi_{(\tilde\beta',R',\tilde\iota_{R'})}\,\colon R'\to  Q'$ fulfilling 
$$\tilde\alpha'\circ \tilde\varphi_{(\beta',R',\iota_{R'})} = 
\tilde\varphi\circ \tilde\beta'\qq 
\tilde\iota_{Q'}\circ \tilde\varphi_{(\tilde\beta',R',\tilde\iota_{R'})} 
=\tilde \iota_{R'}
\eqno £4.9.2;$$
actually, according to conditions~£2.2.1 and~£2.2.2, the right-hand equality is equivalent to say that $ \tilde\varphi_{(\tilde\beta',R',\tilde\iota_{R'})}$ is an 
$\widetilde\A\-$morphism and, in this case, this $\widetilde\A\-$morphism is uniquely determined by the left-hand equality. Denote by~$\widetilde\G$ the
 subcategory of $\widetilde\B$ image of $\G$ and assume that
\smallskip
\noindent
£4.9.3\quad {\it If $\tilde\varphi$ is a $\widetilde\G\-$isomorphism then 
$\tilde\varphi_{(\tilde\beta',R',\tilde\iota_{R'})}$ is also a 
$\widetilde\G\-$isomorphism\/}.

\medskip
£4.10. At this point, we actually get a functor $\tilde I\,\colon \widetilde\B\to \Set$  mapping any $\widetilde\B\-$object $Q$ on  the finite set $\tilde I_Q$ and any $\widetilde\B\-$morphism $\tilde\varphi\,\colon R\to Q$ on the map $\tilde I_{\tilde\varphi}\,\colon \tilde I_R\to \tilde I_Q$ defined above. Then, in the {\it semidirect 
product\/} $\tilde I\!\rtimes \tilde\B$ of $\tilde\B$ by $\tilde I$  (cf.~£2.5)
we simply denote by $(Q',\tilde\alpha',Q)$ the 
$\tilde I\!\rtimes \tilde\B\-$object formed by $Q$ and by 
$(\tilde\alpha',Q',\tilde\iota')\in \tilde I_Q\,,$ whereas a 
$\tilde I\!\rtimes \tilde\B\-$morphism is a pair
$$(\tilde\varphi',\tilde\varphi) : (R',\tilde\beta',R)\too (Q',\tilde\alpha',Q)
\eqno £4.10.1\phantom{.}$$
formed by a $\tilde\B\-$morphism $\tilde\varphi\,\colon R\to Q$ and by
an $\tilde\A\-$morphism $\tilde\varphi'\,\colon R'\to Q'$ fulfilling 
$\tilde\alpha'\circ \tilde\varphi' = \tilde\varphi\circ \tilde\beta'\,.$

\medskip
£4.11.  In this section, as an application of Theorem~£3.5 above we  prove that,  for any {\it contravariant\/} functor 
 $\tilde\frak a\,\colon\widetilde\B\to \Ab\,,$ the {\it standard complex\/} 
 with $n\-$term $\Bbb C^n_{\widetilde\G} (\widetilde\B, \tilde\frak a^{\ad}\circ \widetilde{\frak m}_P)$ (cf.~£2.3),  and with the usual {\it differential map\/}, (cf.~£2.3.2) is {\it homotopically trivial\/}; for this purpose, we consider the 
{\it quintuple\/}~$\H$ formed by
\smallskip
\noindent
£4.11.1\quad {\it The  trivial bi-interior structure over $\widetilde\B\,.$\/}
\smallskip
\noindent
£4.11.2\quad {\it The functor  $\tilde I\,\colon \widetilde\B\too \Set\,.$\/}
\smallskip
\noindent
£4.11.3\quad {\it The functor $\frak n\,\colon \tilde I\!\rtimes \widetilde\B \too \widetilde\A\i \widetilde\B$ sending any $\tilde I\!\rtimes \widetilde\B \-$object 
$(Q',\tilde\alpha',Q)$ to $Q'$ and any 
$\tilde I\!\rtimes \widetilde\B \-$morphism $(\tilde\varphi',\tilde\varphi)\,\colon (R',\tilde\beta',R)\too (Q',\tilde\alpha',Q)$ to the $\widetilde\A\-$morphism 
$\tilde\varphi'\,.$\/}
\smallskip
\noindent
£4.11.4\quad {\it The natural map $\nu\,\colon \frak n\too \id_{\widetilde\B}$ sending any $\tilde I\!\rtimes \widetilde\B \-$object $(Q',\tilde\alpha',Q)$ to the 
$\widetilde\B \-$morphism
$\tilde\alpha'\,\colon Q'\too Q\,.$\/}
\eject
\smallskip
\noindent
This quintuple is clearly a {\it homotopic system\/} for the triple 
$(\widetilde\A,\widetilde\B,\widetilde\G)$ (cf.~£2.6).

\bigskip
\noindent
{\bf Proposition~£4.12.} {\it With the notation above, assume that conditions~{\rm£4.8.1} and~{\rm £4.9.3} are fulfilled. Then, for any contravariant functor
$\tilde\frak a\,\colon \tilde\B\to \Ab\,,$ the contravariant functor 
$\tilde\frak a^{\ad}\circ \widetilde{\frak m}_P$ is $\H\-$split. 
In particular, for any $n\ge 1$ we have
$$\Bbb H^n_{\tilde\G} (\tilde\B,\tilde\frak a^{\ad}\circ \widetilde{\frak m}_P) = \{0\}
\eqno £4.12.1.$$\/}

\par
\noindent
{\bf Proof:} According to definition of {\it $\H\-$split\/} in~£2.10, we have to exhibit a {\it na-tural section\/} of the {\it natural map\/} 
$\Delta_\H (\tilde\frak a^{\ad}\circ \widetilde{\frak m}_P)$ (cf.~£2.10.1); explicitly, for any $\widetilde\B\-$object~$Q$ and any element $a = 
\sum_{(\tilde\alpha',Q',\tilde\iota_{Q'})\in \tilde I_Q} a_{(\tilde\alpha',Q',\tilde\iota_{Q'})}$ in
$$(\tilde\frak a^{\ad}\circ \widetilde{\frak m}_P) (Q) = 
\prod_{(\tilde\alpha',Q',\tilde\iota_{Q'})\in \tilde I_Q} \tilde\frak a (Q')
\eqno £4.12.2,$$
we have (cf.~£2.10.3)
$$\eqalign{\Delta_\H &(\tilde\frak a^{\ad}\circ \widetilde{\frak m}_P)_Q (a) = 
\sum_{(\tilde\alpha',Q',\tilde\iota_{Q'})\in \tilde I_Q} \big((\tilde\frak a^{\ad}\circ \widetilde{\frak m}_P) * \nu\big)_{(Q',\tilde\alpha',Q)} (a)\cr
&= \sum_{(\tilde\alpha',Q',\tilde\iota_{Q'})\in \tilde I_Q} \big((\tilde\frak a^{\ad}\circ \widetilde{\frak m}_P)(\tilde\alpha')\big) (a)\cr
&= \sum_{(\tilde\alpha',Q',\tilde\iota_{Q'})\in \tilde I_Q}\,
\sum_{(\tilde\alpha'',Q'',\tilde\iota_{Q''})\in \tilde I_{Q'}} 
\big(\tilde\frak a (\tilde\alpha'''')\big) (a_{(\tilde\alpha''',Q''',\tilde\iota_{Q'})})\cr}
\eqno £4.12.3\phantom{.}$$ 
where we set $\tilde I_{\tilde\alpha'} (\tilde\alpha'',Q'',\tilde\iota_{Q''}) = 
(\tilde\alpha''',Q''',\tilde\iota_{Q'''})$ for any $(\tilde\alpha'',Q'',\tilde\iota_{Q''})\in \tilde I_{Q'}$ and $\tilde\alpha''''\,\colon Q''\to Q'''$ is the 
$\widetilde\A\-$morphism fulfilling $\tilde\alpha'''\circ \tilde\alpha'''' 
= \tilde\alpha'\circ \tilde\alpha''$ (cf.~£4.9.2).

\smallskip
Now we claim that the group homomorphism
$$\theta_Q : \big(\H (\tilde\frak a^{\ad}\circ \widetilde{\frak m}_P)\big)(Q) \too (\tilde\frak a^{\ad}\circ \widetilde{\frak m}_P)(Q)
\eqno £4.12.4\phantom{.}$$
mapping any element in $\big(\H (\tilde\frak a^{\ad}\circ \widetilde{\frak m}_P)\big)(Q) = {\displaystyle \prod_{(\tilde\alpha',Q',\tilde\iota_{Q'})\in \tilde I_Q}} (\tilde\frak a^{\ad}\circ \widetilde{\frak m}_P)(Q')$
$$b =\sum_{(\tilde\alpha',Q',\tilde\iota_{Q'})\in \tilde I_Q}\,
\sum_{(\tilde\alpha'',Q'',\tilde\iota_{Q''})\in \tilde I_{Q'}} 
b_{(\tilde\alpha'',Q'',\tilde\iota_{Q''})}
\eqno £4.12.5\phantom{.}$$
on ${\displaystyle \sum_{(\tilde\alpha',Q',\tilde\iota_{Q'})\in \tilde I_Q}} 
b_{(\widetilde{\rm id}_{Q'},Q',\tilde\iota_{Q'})}$ is a section of 
$\Delta_\H (\tilde\frak a^{\ad}\circ \widetilde{\frak m}_P)_Q\,,$ and that the cor-respondence sending $Q$ to $\theta_Q$ is {\it natural\/}. Indeed, if in~£4.12.3
above we set $c_{(\tilde\alpha'',Q'',\tilde\iota_{Q''})} = \big(\tilde\frak a (\tilde\alpha'''')\big) (a_{(\tilde\alpha''',Q''',\tilde\iota_{Q'})})$ for any 
$(\tilde\alpha'',Q'',\tilde\iota_{Q''})\in \tilde I_{Q'}\,,$
we get
$$\theta_Q\big(\Delta_\H (\tilde\frak a^{\ad}\circ \widetilde{\frak m}_P)_Q (a)\big) = {\displaystyle \sum_{(\tilde\alpha',Q',\tilde\iota_{Q'})\in \tilde I_Q}} c_{(\widetilde{\rm id}_{Q'},Q',\tilde\iota_{Q'})} = a
\eqno £4.12.6\phantom{.}$$
since the equalities $Q'' = Q'$ and $\tilde\alpha'' =  \widetilde{\rm id}_{Q'}$ force 
$Q''' = Q'\,,$ $\tilde\alpha''' = \tilde\alpha'$ and~$\tilde\alpha'''' = \widetilde{\rm id}_{Q'}\,.$`
\eject

\smallskip
Moreover, on the one hand, for any $\widetilde\B\-$morphism $\tilde\varphi\,
\colon R\to Q\,,$ we have
$$\eqalign{\big((\tilde\frak a^{\ad}\circ \widetilde{\frak m}_P)&(\tilde\varphi) 
\circ \theta_Q\big) (b)\cr 
&= \big((\tilde\frak a^{\ad}\circ \widetilde{\frak m}_P)(\tilde\varphi) \big)
\big(\sum_{(\tilde\alpha',Q',\tilde\iota_{Q'})\in \tilde I_Q} 
b_{(\widetilde{\rm id}_{Q'},Q',\tilde\iota_{Q'})}\big)\cr
&= \sum_{(\tilde\beta',R',\tilde\iota_{R'})\in \tilde I_R} 
\big(\tilde\frak a (\tilde\varphi')\big)
(b_{(\widetilde{\rm id}_{Q'},Q',\tilde\iota_{Q'})})\cr}
\eqno £4.12.7\phantom{.}$$
where if $\tilde I_{\tilde\varphi} (\tilde\beta',R',\tilde\iota_{R'}) 
= (\tilde\alpha',Q',\tilde\iota_{Q'})$ then $\tilde\varphi'\,\colon R'\to Q'$ is the
$\widetilde\A\-$morphism fulfilling $\tilde\alpha'\circ \tilde\varphi' 
= \tilde\varphi\circ \tilde\beta'$ (cf.~£4.9.2); on the other hand, we get
$$\eqalign{\Big(\big(\H (\tilde\frak a^{\ad}&\circ \widetilde{\frak m}_P)\big)
(\tilde\varphi)\Big) (b)\cr
& = \sum_{(\tilde\beta',R',\tilde\iota_{R'})\in \tilde I_R}\,
\sum_{(\tilde\beta'',R'',\tilde\iota_{R''})\in \tilde I_{R'}} 
\big(\tilde\frak a (\tilde\varphi'')\big) (b_{(\tilde\alpha'',Q'',\tilde\iota_{Q''})})\cr}
\eqno £4.12.8\phantom{.}$$
where if $\tilde I_{\tilde\varphi} (\tilde\beta',R',\tilde\iota_{R'}) 
= (\tilde\alpha',Q',\tilde\iota_{Q'})\,,$  $\tilde\varphi'\,\colon R'\to Q'$ is the
$\widetilde\A\-$morphism fulfilling $\tilde\alpha'\circ \tilde\varphi' 
= \tilde\varphi\circ \tilde\beta'\,,$ and $\tilde I_{\tilde\varphi}' 
(\tilde\beta'',R'',\tilde\iota_{R''})  = (\tilde\alpha'',Q'',\tilde\iota_{Q''})\,,$
 then $\tilde\varphi''\,\colon R''\to Q''$ is the $\widetilde\A\-$morphism fulfilling 
 $\tilde\alpha''\circ \tilde\varphi'' = \tilde\varphi'\circ \tilde\beta''$
(cf.~£4.9.2); consequently, for any $(\tilde\beta',R',\tilde\iota_{R'})
\in \tilde I_R\,,$ the term labeled by $(\widetilde{\rm id}_{R'},R',\tilde\iota_{R'})$ in the second sum coincides with $\big(\tilde\frak a (\tilde\varphi')\big)
(b_{(\widetilde{\rm id}_{Q'},Q',\tilde\iota_{Q'})})$ since the equalities $R'' = R'$ 
and $\tilde\beta'' = \widetilde{\rm id}_{R'}\,,$ and the fact that $\tilde\varphi'$
is already an $\widetilde\A\-$morphiusm force $Q'' = Q'\,,$ 
$\tilde\alpha'' = \widetilde{\rm id}_{Q'}$ and $\tilde\varphi'' = \tilde\varphi'\,.$
Finally, we obtain
$$\Big(\theta_R\circ \big(\H (\tilde\frak a^{\ad}\circ \widetilde{\frak m}_P)\big)
(\tilde\varphi)\Big) (b) = \big((\tilde\frak a^{\ad}\circ \widetilde{\frak m}_P)
(\tilde\varphi)  \circ \theta_Q\big) (b)
\eqno £4.12.9.$$
We are done.

\medskip
£4.13. Let us give some explanation on the application of this theorem in 
the proof of [4,~Theorem~6.22]. In this case, as in~£4.2 above, we are considering
a finite $p\-$group $P\,,$ a Frobenius $P\-$category $\F$ and a nonempty set
$\frak X$ of $\F\-$selfcentralizing subgroups of $P$  which contains any subgroup of $P$ admitting an $\F\-$morphism from some subgroup in~$\frak X\,.$ Then, $\B$
is the {\it full\/} subcategory $\F^{^\frak X}$ of $\F$ over $\frak X\,,$ endowed with
the usual {\it interior structure\/}, $\A$ is $\F_{\!P}^{^\frak X}$ and $\G$ coincides with
the {\it interior structure\/}; we already know that $\F^{^\frak X}$ is a {\it multiplicative
$\F_{\!P}^{^\frak X}\-$category\/} and condition~£4.9.3 is clearly fulfilled. Thus,
considering the {\it center  contravariant\/} functor $\tilde\frak z^{_\frak X}\,\colon 
\widetilde\F^{^\frak X}\to \Ab$ mapping $Q\in \frak X$ on $Z(Q)\,,$ for any $n\ge 1$ we get
$$\Bbb H^n (\widetilde\F^{^\frak X}, (\tilde\frak z^{_\frak X})^{\ad}\circ \widetilde{\frak m}_P)
 = \{0\}
\eqno £4.13.1.$$

\medskip
£4.14. Our argument in  [4,~Theorem~6.22] goes by induction on $\vert \frak X\vert\,,$
so that we consider a minimal element $U$ in $\frak X$ {\it fully normalized\/} in $\F$
[3,~2.6]; then, denoting by
$\tilde\frak z_U^{_\frak X}\,\colon \widetilde\F^{^\frak X}\to \Ab$ the functor sending
$Q\in \frak X$ to $Z(Q)$ or to $\{0\}$ according to $Q$ is $\F\-$isomorphic to $U$ or not,
we have an obvious ``inclusion'' $\tilde\frak z_U^{_\frak X}\i 
(\tilde\frak z^{_\frak X})^{\ad}\circ \widetilde{\frak m}_P$ and are reduced to prove that some {\it $\big((\tilde\frak z^{_\frak X})^{\ad}\circ \widetilde{\frak m}_P\big)\big/\tilde\frak z_U^{_\frak X}\-$valued $2\-$cocycle\/} over $\widetilde\F^{^\frak X}$ is a {\it $2\-$boundary\/}. But, from equality~£4.13.1 we obtain
$$\Bbb H^2 \big(\widetilde\F^{^\frak X}\!, \big((\tilde\frak z^{_\frak X})^{\ad}\circ \widetilde{\frak m}_P\big)\big/\tilde\frak z_U^{_\frak X}\big)
\cong \Bbb H^3 \big(\widetilde\F^{^\frak X}\!, \tilde\frak z_U^{_\frak X}\big)
\eqno £4.14.1$$
and it follows from [1, Proposition~3.2] that, setting $\N = N_\F (U)$ and denoting by
$\frak N$ the set of $Q\in \frak X$ contained in $N_P (U)$ and by $\tilde\frak z^{_ {\frak N}}_U$ the restriction of $\tilde\frak z^{_ {\frak X}}_U$ to $ \widetilde\N^{^\frak N}\,,$
the ``inclusion'' $\widetilde\N^{^\frak N}\i \widetilde\F^{^\frak X}$ induces [4,~Proposition~6.9]
$$\Bbb H^3 (\tilde\F^{^{\frak X}},\tilde\frak z^{_ {\frak X}}_U)\cong 
 \Bbb H^3 \big(\tilde\N^{^{\frak N}}\!,\tilde\frak z^{_ {\frak N}}_U\big)
 \eqno £4.14.2.$$
Finally, {\it mutatis mutandis\/}  from equality~£4.13.1 we also obtain
$$\Bbb H^2 \big(\widetilde\N^{^\frak N}\!, \big((\tilde\frak z^{_\frak N})^{\ad}\circ \widetilde{\frak m}_P\big)\big/\tilde\frak z_U^{_\frak N}\big)
\cong \Bbb H^3 \big(\widetilde\N^{^\frak N}\!, \tilde\frak z_U^{_\frak N}\big)
\eqno £4.15.1\phantom{.}$$
and succeed prove that our {\it $2\-$cocycle\/} can be translated to the analogous
{\it $2\-$cocycle\/} for $\widetilde\N^{^\frak N}\,;$ here, the existence of the {\it $\F\-$localizer\/} $L_\F (U)$ [3,~Theorem~18.6] completes our proof [4,~Proposition~6.19].

 \bigskip
\noindent
{\bf £5. The compatible complement case\/}
\medskip
£5.1. Let $P$ be a finite $p\-$group and $\F$ a Frobenius $P\-$category [3,~2.8]; we are interested in a family  of commutative square diagrams in the {\it additive cover\/} 
$\ad(\widetilde\F)$ (cf.~£4.1) of the {\it exterior quotient\/}~$\widetilde\F$ (cf.~£4.2) --- called {\it special $\ad(\widetilde\F)\-$squares\/}. Explicitly, for any triple of subgroups $Q\,,$ $R$ and $T$ of $P\,,$ and any pair of $\widetilde\F\-$morphisms 
$\tilde\alpha\,\colon R\to Q$ and $\tilde\beta\,\colon T\to Q\,,$ we choose  a pair of representatives   $\alpha$ of $\tilde\alpha$ and $\beta$ of $\tilde\beta\,,$ and  a set of representatives $W\i Q$ for the set of double classes $\alpha(R)\backslash Q/\beta(T)\,;$ moreover, for any $w\in W\,,$ we set 
 $U_w = \alpha(R)^w\cap \beta(T)$ and respectively denote by
$$\alpha_w : U_w\too R\qq \beta_w : U_w\too T
\eqno £5.1.1\phantom{.}$$
the $\F\-$morphisms mapping $u\in U_w$ on $\alpha^*(wuw^{-1})$ and on 
$\beta^*(u)\,,$ where $\alpha^*\,\colon \alpha(R)\cong R$ and 
$\beta^*\,\colon \beta(T)\cong T$ are the inverse maps of the respective isomorphisms induced by $\alpha$ and $\beta\,.$

\medskip
£5.2. Consider the  $\ad(\widetilde\F)\-$object $U = \bigoplus_{w\in W} U_w$
and the two $\ad(\widetilde\F)\-$mor-phisms $\tilde\gamma\,\colon U\to R$ and  
$\tilde\delta \,\colon U\to T$ respectively determined by the families 
$\{\tilde\alpha_w\}_{w\in W}$ and  $\{\tilde\beta_w\}_{w\in W}$ (cf.~£4.1);
then, it is quite clear that the following\break
\eject
\noindent 
$\ad(\widetilde\F)\-$diagram is commutative
$$\matrix{&&Q\cr
&{\tilde\alpha\atop}\hskip-5pt\nearrow\hskip-5pt&&\hskip-5pt\nwarrow\hskip-5pt{\tilde\beta\atop}\cr
R\hskip-5pt&&&&\hskip-5pt T\cr
&{\atop \tilde\gamma}\hskip-5pt\nwarrow\hskip-5pt&&\hskip-5pt\nearrow\hskip-5pt{\atop \tilde\delta}\cr
&&U\cr}
\eqno £5.2.2.$$
Let us call {\it special $\ad(\widetilde\F)\-$square\/} any {\it direct sum\/} in 
$\ad(\widetilde\F)$ of such 
commutative $\ad(\widetilde\F)\-$diagrams. 

\medskip
£5.3.  Let $\frak a\,\colon \widetilde\F\to \O\-\mod$ be a {\it contravariant\/}
functor and consider the  {\it additive contravariant\/} functor $\frak a^{_{\ad}}\,\colon \ad(\widetilde\F)\to \O\-\mod$ (cf.~£4.8.4); we say that  a functor 
$\frak a^\circ\,\colon \widetilde\F\to \O\-\mod$ is a {\it compatible complement\/} of $\frak a$
  if it fulfills the following three conditions:
    \smallskip
  \noindent
£5.3.1\quad  {\it   For any subgroup $Q$ of $P$ we have 
$\frak a^\circ (Q) = \frak a (Q)\,.$\/}
  \smallskip
  \noindent
£5.3.2\quad {\it  For any $\widetilde\F\-$morphism $\tilde\varphi\,\colon R\to Q$ we have $\frak a^\circ (\tilde\varphi)\circ\frak a (\tilde\varphi) = 
 {\vert Q\vert\over\vert R\vert}\.{\rm id}_{\frak a (Q)}\,.$\/}
   \par
  \noindent
£5.3.3\quad  {\it For any special $\ad(\widetilde\F)\-$square
$\matrix{&&Q\cr
&{\tilde\varphi\atop}\hskip-5pt\nearrow\hskip-5pt&&\hskip-5pt\nwarrow\hskip-5pt{\tilde\psi\atop}\cr
R\hskip-10pt&&&&\hskip-10pt T\cr
&{\atop\tilde\gamma}\hskip-5pt\nwarrow\hskip-5pt&&\hskip-5pt\nearrow\hskip-5pt{\atop\tilde\delta}\cr
&&U\cr}$
we have the commutative $\O\-\mod\-$square
$\matrix{&&\frak a^{_{\ad}} (Q)\cr
&{(\frak a^\circ)^{_{\ad}} (\tilde\varphi)\atop}\hskip-5pt\nearrow\hskip-1pt&
&\hskip-5pt\searrow\hskip-5pt{\frak a^{_{\ad}}(\tilde\psi)\atop}\cr
\frak a^{_{\ad}} (R)\hskip-40pt&&&&\hskip-40pt \frak a^{_{\ad}} (T)\cr
&{\atop\frak a^{_{\ad}} (\tilde\gamma)}\hskip-5pt\searrow\hskip-10pt&
&\hskip-1pt\nearrow\hskip-5pt{\atop(\frak a^\circ)^{_{\ad}} (\tilde\delta)}\cr
&&\frak a^{_{\ad}} (U)\cr}
\,.$ \/}
  \par
  \noindent
Note that $\frak a^\circ (\tilde\varphi) = \frak a (\tilde\varphi^{-1})$ for any 
$\widetilde\F\-$isomorphism $\tilde\varphi\,.$

\medskip
£5.4. Actually, the {\it special $\ad(\widetilde\F)\-$squares\/} are the {\it image\/} of the 
{\it pull-backs\/} of a suitable extension of the category $\ad(\widetilde\F)\,;$ in order to
explicit this extension we need some notation.  Recall that any Frobenius  $P\-$category  comes from a {\it basic $P\times P\-$set\/} [3,~Propositions~21.9 and~21.12] and let $\Omega$ be  an {\it $\F\-$basic $P\times P\-$set\/} [3,~21.3]; as in [3,~22.2], let us denote by~$G$ the group of $\{1\}\times P\-$set automorphisms of  ${\rm Res}_{\{1\}\times P}(\Omega)$ and identify $P$ with the image  of~$P\times \{1\}$ in~$G\,.$ Then, 
we consider the categories $\T$ and $\T_P$ where the objects still are the subgroups 
of~$P$ and the morphisms are given by the respective {\it transporters\/} in 
$G\times P$ and in $P\times P\,;$ explicitly, for any pair of subgroups 
 $Q$ and $R$ of $P\,,$ we have
 $$\T (Q,R) = T_G (R,Q)\times P\qq \T_P (Q,R) 
 = T_P (R,Q)\times P
 \eqno £5.4.1,$$
 where $T_G$ and $T_P$ respectively denote the {\it transporters\/} in $G$ and in $P\,,$ the composition being defined by the products in $G$ and in $P\,;$  it is  clear that $\T$ and $\T_P$ are ordered categories and that $\T$ is a 
{\it $\T_P\-$category\/}.
 
 \medskip
£5.5. Thus, by the very definition of {\it $\F\-$basic $P\times P\-$sets\/} [3,~21.3], $\F$ and~$\F_{\!P}$ are actually respective {\it quotients\/} of  $\T$ and $\T_P\,;$ more precisely, $\T$~admits the {\it co-interior\/} structure which  maps any subgroup $Q$ 
of $P$ on $C_G(Q)\times P$ and then $\F$ is just the corresponding 
{\it exterior quotient\/}. Similarly, $\T$ admits the {\it interior\/} structure $\I$ mapping 
 $Q$  on $Q\times P\,;$ let us denote by~$\widetilde\T^{^\I}$ the corresponding {\it exterior quotient\/} and by $\tilde\frak e\,\colon \widetilde\T^{^\I}\to 
\widetilde\F$ the obvious canonical functor. Finally, $\T$~also admits a {\it bi-interior\/}
structure mapping $Q$ on $\I (Q) = Q\times P$ and on~$\I^\circ (Q) =
C_G(Q)\times \{1\}\,;$ note that we have $\widetilde\T = \widetilde\F\,.$

\bigskip
\noindent
{\bf Proposition~£5.6.} {\it  The additive cover $\ad\big(\widetilde\T^{^\I}\big)$ of 
$\widetilde\T^{^\I}$ admits pull-backs and the  special $\ad(\widetilde\F)\-$squares are  
the image by $\tilde\frak e$ of all the  
$\ad\big(\widetilde\T^{^\I}\big)\-$pull-backs\/}.

\medskip
\noindent
{\bf Proof:} In order to prove that  $\ad\big(\widetilde\T^{^\I}\big)$ admits {\it pull-backs\/}, 
it suffices to consider any triple of subgroups $Q\,,$ $R$ and $T$ of $P\,,$ and any pair of 
$\widetilde\T^{^\I}\-$morphisms  $\tilde x\,\colon R\to Q$ and $\tilde y\,\colon T\to Q\,.$  
In this case, choose  a pair of representatives   $x\in T_G (R,Q)$ of $\tilde x$ and $y\in T_G (T,Q)$ of $\tilde y\,,$ and  a set of representatives $W\i Q$ for the set of double classes ${}^x R\backslash Q/{}^y T\,;$ moreover, 
for any $w\in W\,,$  set  $U_w = {}^x R\cap{}^{wy}(T)\,,$  respectively denote by
$$x_w : U_w\too R\qq y_w : U_w\too T
\eqno £5.6.1\phantom{.}$$
the $\T\-$morphisms mapping $u\in U_w$ on  $u^x$ and on $u^{wy}\,,$ and
consider the  two $\ad\big(\widetilde\T^{^\I}\big)\-$morphisms
$$\bigoplus_{w\in W}\tilde x_w : \bigoplus_{w\in W} U_w\too R\qq \bigoplus_{w\in W}\tilde y_w : \bigoplus_{w\in W} U_w\too T
\eqno £5.6.2\phantom{.}$$
respectively determined by the families $\{\tilde x_w\}_{w\in W}$ and 
$\{\tilde y_w\}_{w\in W}\,.$ Then, we claim that the following 
$\ad\big(\widetilde\T^{^\I}\big)\-$diagram is a {\it pull-back\/}
$$\matrix{&&Q\cr
&{\tilde x\atop}\hskip-5pt\nearrow\hskip-50pt&&\hskip-50pt\nwarrow\hskip-5pt{\tilde y\atop}\cr
R\hskip-70pt&&\phantom{\Big\uparrow}&&\hskip-45pt T\cr
&{\atop { \bigoplus_{w\in W}}\tilde x_w}\hskip-5pt\nwarrow\hskip-10pt&&\hskip-10pt\nearrow\hskip-2pt{\atop { \bigoplus_{w\in W}}\tilde y_w}\cr
&&{\bigoplus_{w\in W}} U_w\cr}
\eqno £5.6.3.$$

\smallskip
Indeed,  if $S$ is a subgroup of $P$ and we have two 
$\widetilde\T^{^\I}\!\-$morphisms 
$\tilde a\,\colon S\to R$ and $\tilde b\,\colon  S\to T$ fulfilling $\tilde x\.\tilde a = 
\tilde y\.\tilde b$ then, choosing representatives $a$ in $\tilde a$ and 
$b$ in $\tilde b\,,$ there exists $u\in Q$ such that $xa = uyb$ and it is easily checked that the double class of $u$ in 
${}^x R\backslash Q/{}^y T $ does not depend on our choices of $a$ and~$b\,;$ moreover, denoting by $w\in W$ the representative of the double class of $u$ and choosing $r\in R$
and $t\in T$ fulfilling $u = ({}^x r)^{-1} w({}^y t)\,,$ we get
$$xra = wytb
\eqno £5.6.4\phantom{.}$$
and therefore ${}^{xra } S$ is contained in $U_w\,;$ thus, we obtain an 
$\widetilde\T^{^\I}\!\-$morphism $\widetilde{xra}\,\colon S\to U_w$ and therefore an 
$\ad(\widetilde\T^{^\I})\-$morphism
$$S\too {\displaystyle\bigoplus_{w'\in W}} U_{w'}
\eqno £5.6.5\phantom{.}$$
determined by the element  $w\in W$  and by the $\widetilde\T^{^\I}\!\-$morphism 
$\widetilde{xra}\,\colon S\to U_w$ which clearly fulfills
$$\tilde x_w\circ\widetilde{xra} = \tilde a\qq  \tilde y_w\circ\widetilde{xra} = \tilde b
\eqno £5.6.6.$$

\smallskip
Conversely, for any $\ad(\widetilde\T^{^\I})\-$morphism 
$S\to {\displaystyle\bigoplus_{w''\in W}} U_{w''}$ determined by a suitable $w'\in W$ and by a $\widetilde\T^{^\I}\!\-$morphism 
$\tilde c\,\colon S\to U_{w'}$ fulfilling  $\tilde x_{w'}\circ\tilde c = \tilde a$ and  
$\tilde y_{w'}\circ\tilde c = \tilde b,$ it is easily checked that $w' = w$ and that 
$\tilde c = \widetilde{xra}\,.$ Now, the last statement is clear.

\medskip
\noindent
{\bf Remark £5.7.} According to this proposition, to say that a {\it contravariant\/} functor 
$\frak a\,\colon \widetilde\F\to \O\-\mod$ admits a {\it compatible complement\/}
$\frak a^\circ\,\colon \widetilde\F\to \O\-\mod\,,$ in terms of~[2] is equivalent to say that
the functor
$$(\frak a^\circ\circ\tilde\frak e)^{^{\ad}} : \ad\big(\widetilde\T^{^\I}\big)\too \O\-\mod 
\eqno £5.7.1\phantom{.}$$ 
is a {\it  cohomological Mackey complement\/} of the {\it contravariant\/}
functor 
$$(\frak a\circ\tilde\frak e)^{^{\ad}} : \ad\big(\widetilde\T^{^\I}\big)\too \O\-\mod 
\eqno £5.7.2.$$

 \medskip
 £5.8.  On the other hand, we consider a subcategory $\G$ of $\T$
having the same objects and only $\T\-$isomorphisms; we assume that $\G$ contains all the $\T_P\-$isomorphisms and fulfills
$$\I(Q)\;\I^\circ (Q) = Q\.C_G (Q)\times P\i \G (Q)
\eqno £5.8.1\phantom{.}$$
for any subgroup $Q$ of $P\,.$ Recall that, for any  subgroup $Q$ of $P\,,$ the stabilizer $(Q\times P)_\omega$ of $\omega\in \Omega$ in $Q\times P$ has the form [3,~21.3]
$$\Delta_{t^Q_\omega}(Q_\omega) = 
\{\big({}^{t^Q_\omega} v,v\big)\}_{v\in Q_\omega}
\eqno £5.8.2\phantom{.}$$
for some subgroup $Q_\omega$ of $P$ and some $t^Q_\omega \in 
T_G (Q_\omega,Q)\,;$ note that, for any $\T\-$morphism 
$(x,u)\,\colon R\to Q$ and any $\omega\in \Omega\,,$ we have
$${}^u R_\omega \i Q_{x\.\omega\.u^{-1}}\qq x\. t^R_\omega \equiv t^Q_{x\.\omega\.u^{-1}}\.u \bmod{C_G(R_\omega)}
\eqno £5.8.3.$$

\medskip
£5.9. As in~£2.6, we consider the following quintuple 
$\H = (\I,\I^\circ,\frak s,\frak n,\nu)$  
\smallskip
\noindent
£5.9.1\quad {\it  $\I(Q) = Q\times P$ and $\I^\circ (Q) = C_G(Q)\times\{1\}$ for any subgroup $Q$ of~$P\,.$\/}
\smallskip
\noindent
£5.9.2\quad {\it The functor $\frak s\,\colon \T\to \Set$ maps any  subgroup $Q$ of~$P$ on the set $\Omega$ and any $\T\-$morphism $(x,u)\,\colon R\to Q$ on the  map sending $\omega\in \Omega$
 to $x\.\omega\.u^{-1}\,.$\/}
\smallskip
\noindent
£5.9.3\quad {\it The functor $\frak n\,\colon \frak s\!\rtimes \T\to \widetilde\T_P =  
\widetilde\F_P$ maps any  $\frak s\!\rtimes \T\-$object $(\omega,R)$ on $R_\omega\i P$ and any $\frak s\!\rtimes \T\-$morphism 
$$(x,u) : (\omega,R)\too (x\.\omega\.u^{-1},Q)$$
on the exterior class of the conjugation by $u\in P$ from $R_\omega$ to $Q_{x\.\omega\.u^{-1}}\,.$\/}
\eject
\smallskip
\noindent
£5.9.4\quad {\it The natural map $\nu\,\colon \frak n\to \tilde\frak p$  sends  
any $\frak s\!\rtimes \T\-$object
$(\omega,Q)$ to  the bi-exterior class of $(t^Q_\omega,1)\,\colon Q_\omega\to Q\,.$\/}
\smallskip
\noindent
This quintuple is a {\it homotopic system\/} for the triple $(\T_P,\T,\G)\,;$ indeed, the functoriality of $\frak s$ comes from the very definition of $\T\,;$ since $\widetilde\G$
contains all the $\widetilde\T_P\-$isomorphisms, $\frak n$ maps $\frak s_\G\rtimes \G$
in $\widetilde\G\,;$ the {\it functoriality\/} of $\frak n$ and the {\it naturality\/} of $\nu$ follow from~£5.8.3.

\bigskip
\noindent
{\bf Theorem~£5.10.} {\it With  the notation and the hypothesis above,  any contra-variant functor $\frak a\,\colon \widetilde\F\to \O\-\mod$  admitting a compatible complement $\frak a^\circ$ is $\H\-$split. In particular, for any $n\ge 1$ we have  $\Bbb H_{\tilde\G}^n (\widetilde\F,\frak a) = \{0\}\,.$\/}

\medskip
\noindent
{\bf Proof:} Setting $\hat\frak a = \frak a\circ \frak e_{\I,\I^\circ}\,,$ it follows from~£2.10 that there is a {\it natural map\/} $\Delta_\H (\frak a)\,\colon
\hat\frak a \to \H(\frak a)$
sending any subgroup $Q$ of $P$ to the group homomorphism
$$\Delta_\H (\frak a)_Q : \frak a (Q) \too \big(\prod_{\omega\in \Omega} \frak a(Q_\omega)\big)^{Q\times P}
\eqno £5.10.1\phantom{.}$$
mapping $a\in \frak a (Q)$ on $\big(\Delta_\H (\frak a)_Q \big) (a) = \sum_{\omega\in \Omega} \big(\frak a (\,\widetilde{ t^Q_\omega}\,)\big)(a)\,.$

\smallskip
First of all note that, for any subgroup $Q$ of $P$ and any $\omega\in \Omega\,,$ we  have (cf~£5.3.1)
$$\big(\frak a^\circ (\,\widetilde{ t^Q_\omega}\,)\big)\Big( \big(\frak a (\,\widetilde{ t^Q_\omega}\,)\big)(a)\Big)
= {\vert Q\vert\over \vert Q_\omega\vert}\.a
\eqno £5.10.2.$$
Then, choosing a set of representatives $\Gamma_Q\i \Omega$ for the set of $Q\times P\-$orbits 
in~$\Omega\,,$ we consider the map
$$\theta_Q : \big(\prod_{\omega\in \Omega} \frak a(Q_\omega)\big)^{Q\times P}\too \frak a (Q)
\eqno £5.10.3\phantom{.}$$
sending any element  $\sum_{\omega\in \Omega} a_\omega$ of $\big(\H (\frak a)\big)(Q)\,,$ where $a_\omega$ belongs to $\frak a(Q_\omega)$ for any $\omega\in \Omega\,,$ to 
$$\theta_Q(\sum_{\omega\in \Omega} a_\omega) = {\vert P\vert\over \vert\Omega\vert}\.\sum_{\omega\in\Gamma_Q}\big(\frak a^\circ (\,\widetilde{ t^Q_\omega}\,)\big)(a_\omega)
\eqno £5.10.4;$$
in particular, for any $a\in \frak a (Q)$ we get (cf.~£5.10.2)
$$\eqalign{\theta_Q\Big(\big(\Delta_\H (\frak a)_Q \big) (a)\Big) &= {\vert P\vert\over \vert\Omega\vert}\.\sum_{\omega\in\Gamma_Q}{\vert Q\vert\over \vert Q_\omega\vert}\.a = {\vert P\vert\over \vert\Omega\vert}\.\sum_{\omega\in\Gamma_Q}{\vert Q\.\omega\. P\vert\over \vert P\vert}\.a = a\cr}
\eqno £5.10.5.$$
Moreover, it easily follows from~£5.8.3 and from condition~£5.3.2 that the map $\theta_Q$ does not depend on the choice of $\Gamma_Q\,.$

\smallskip
Consequently, according to Theorem~£3.5, it suffices to prove that the correspondence sending  any  subgroup $Q$ of $P$ to $\theta_Q$ is a 
{\it natural map\/}. Let $(x,u)\,\colon R\to Q$ be a $\T\-$morphism and denote by $\tilde x\colon R\to Q$ the image of $(x,u)$\break
\eject
\noindent
 in~$\widetilde\F (Q,R)\,;$ we claim that the 
following  diagram is commutative
$$\matrix{\big(\prod_{\omega\in \Omega} \frak a(Q_\omega)\big)^{Q\times P}
&\buildrel \theta_Q\over\too &\frak a (Q)\cr
\hskip-60pt{\scriptstyle (\H(\frak a))  ( x,u)}\big\downarrow&
\phantom{\Big\downarrow}&\big\downarrow{\scriptstyle \frak a (\tilde x)}\hskip-20pt\cr
\big(\prod_{\omega\in \Omega} \frak a(R_\omega)\big)^{R\times P}&\buildrel \theta_R\over\too &\frak a (R)\cr}
\eqno £5.10.6;$$
indeed, it follows from definition~£5.10.4 that
$$(\frak a (\tilde x)\circ\theta_Q)(\sum_{\omega\in \Omega} a_\omega) = 
{\vert P\vert\over \vert\Omega\vert}\.\sum_{\omega\in\Gamma_Q}\big(\frak a (\tilde x)\circ
\frak a^\circ (\,\widetilde{ t^Q_\omega}\,)\big)(a_\omega)
\eqno £5.10.7.$$

\smallskip
But, for any $\omega\in\Gamma_Q\,,$ choosing a set of representatives $V_\omega\i Q$  for the set of double classes 
${}^{ t^Q_\omega}Q_\omega\big\backslash Q\big/\,{}^x R\,,$ so that we have
$$Q = \bigsqcup_{v\in V_\omega} ({}^x R)v^{-1}({}^{ t^Q_\omega}Q_\omega)
\eqno £5.10.8,$$
 and setting $T_{\omega,v}= {}^{vx} R\cap{}^{ t^Q_\omega}Q_\omega$ and  
 $x_{\omega,v} = vx$
 for any $v\in V_\omega\,,$ the pair of $\ad(\widetilde\F)\-$morphisms
$\tilde\gamma_\omega\,\colon\bigoplus_{v\in V_\omega} T_{\omega,v}\to R$
and $\tilde\delta_\omega \,\colon\bigoplus_{v\in V_\omega} T_{\omega,v}\to Q_\omega$
 determined by the families of $\widetilde\F\-$morphisms 
 $$\tilde\gamma_{\omega,v} = \widetilde{ x_{\omega,v}}^{-1} : 
 T_{\omega,v}\too R\qq 
 \tilde\delta_{\omega,v} =\widetilde{ t^Q_\omega}^{-1} : 
 T_{\omega,v}\too Q_\omega
 \eqno £5.10.9\phantom{.}$$
  form a {\it special $\ad(\widetilde\F)\-$square\/} (cf.~£5.1 and~£5.2)
$$\matrix{&&Q\cr
&{\widetilde{ t^Q_\omega}\atop}\hskip-5pt\nearrow\hskip-15pt&&\hskip-15pt\nwarrow\hskip-5pt{\tilde x\atop}\cr
Q_\omega\hskip-20pt&&&&\hskip-15pt R\cr
&{\atop \tilde\delta_\omega}\hskip-5pt\nwarrow\hskip-10pt&
&\hskip-10pt\nearrow\hskip-5pt{\atop \tilde\gamma_\omega}\cr
&&\bigoplus_{v\in V_\omega} T_{\omega,v}\cr}
\eqno £5.10.10.$$
Consequently, we still have the commutative $\O\-\mod\-$square
$$\matrix{&&\frak a (Q)\cr
&{\frak a^\circ (\widetilde{ t^Q_\omega})\atop}\hskip-5pt\nearrow\hskip-15pt&
&\hskip-35pt\searrow\hskip-5pt{\frak a(\tilde x)\atop}\cr
\frak a (Q_\omega)\hskip-30pt&&\phantom{\Big\downarrow}&&\hskip-35pt \frak a (R)\cr
&{\atop\frak a^{_{\ad}} (\tilde\delta_\omega)}\hskip-5pt\searrow\hskip-15pt&\phantom{\big\downarrow}
&\hskip-15pt\nearrow\hskip-5pt{\atop(\frak a^\circ)^{^{\ad}} (\tilde\gamma_\omega)}\cr
&&\prod_{v\in V_\omega}\frak a ( T_{\omega,v})\cr}
\eqno £5.10.11\phantom{.}$$ 
and therefore equality~£5.10.7 becomes
$$(\frak a (\tilde x)\circ\theta_Q)(\sum_{\omega\in \Omega} a_\omega) = 
{\vert P\vert\over \vert\Omega\vert}\.\sum_{\omega\in\Gamma_Q} \sum_{v\in V_\omega}
\big(\frak a^\circ (\widetilde{ x_{\omega,v}}^{-1})\circ
\frak a (\,\widetilde{ t^Q_\omega}^{-1})\big)(a_\omega)
\eqno £5.10.12.$$
\eject

\smallskip
On the other hand,  choosing a set of representatives $\Gamma_{\!R}\i \Omega$ for the set of $R\times P\-$orbits 
in~$\Omega\,,$  it follows again from definition~£5.10.4 and from~£5.8.3 that
$$\Big(\theta_R\circ \big(\H(\frak a)\big)  ( x,u)\Big)
(\sum_{\omega\in \Omega} a_\omega) = 
{\vert P\vert\over \vert\Omega\vert}\.\sum_{\omega\in\Gamma_{\!R}}\big(
\frak a^\circ (\,\widetilde{ t^R_\omega}\,)\circ \frak a (\tilde u)\big)(a_{x\.\omega\.u^{-1}})
\eqno £5.10.13.$$
But, for any $\omega\in \Gamma_{\!Q}\,,$ it is clear that we have
$$Q\.\omega\.P =  \bigsqcup_{v\in V_\omega}( {}^xR) v^{-1}\.\omega\. P =\bigsqcup_{v\in V_\omega} xRx^{-1}v^{-1}\.\omega\. P
\eqno £5.10.14\phantom{.}$$
and, in particular, we can choose $\Gamma_{\!R} =\bigsqcup_{\omega\in \Gamma_{\!Q}} x^{-1}(V_\omega)^{-1}\.\omega\.u$
as a set of representatives; in this case, for any $\omega\in \Gamma_{\!Q}$ and any $v\in V_\omega\,,$ the terms respectively corresponding to the pair $(\omega,v)$ and to $x^{-1}v^{-1}\.\omega\.u\in \Gamma_{\!R}$ in the right-hand members of equalities~£5.10.12 and £5.10.13 coincide with each other, as it can be checked from the following commutative diagram
$$\matrix{\frak a (Q_\omega)&\buildrel \frak a (\,\widetilde{ t^Q_\omega}^{-1})\over{\hbox to 35pt{\rightarrowfill}} &\frak a (T_{\omega,v})&\buildrel \frak a^\circ (\widetilde{ x_{\omega,v}}^{-1})\over{\hbox to 65pt{\rightarrowfill}} & \frak a (R)\cr
\Vert&&{\scriptstyle \frak a (\,\widetilde{ t^Q_\omega u})}\big\downarrow\big\uparrow{\scriptstyle \frak a^\circ (\,\widetilde{v t^Q_{v^{-1}\omega} u})}
\hskip-20pt&\phantom{\Big\uparrow}&\Vert\cr
\frak a(Q_{v^{-1}\.\omega}) &\buildrel \frak a (\tilde u)\over{\hbox to 35pt{\rightarrowfill}} & \frak a(R_{x^{-1}v^{-1}\.\omega\.u})&\buildrel \frak a^\circ (\,\widetilde{ t^R_{x^{-1}v^{-1}\.\omega\.u}}\,)\over{\hbox to 65pt{\rightarrowfill}}  & \frak a (R)\cr}
\eqno £5.10.15\,.$$

\medskip
£5.11. The application of this theorem to the proof of [4,~Theorem~£8.10] is quite direct:
the purpose of this result is to lift, to the {\it basic $\F\-$locality\/}~$\L^{^b}$ [3,~22.1], a functor  from the {\it perfect $\F\-$locality\/} $\P$ [3,~17.13] to a suitable quotient 
$\bar\L^{^b}$ of $\L^{^b}\,.$ Then, the {\it Abelian kernel\/} of the functor
from $\L^{^b}$ to $\bar\L^{^b}$ can be {\it filtered\/} in such a way that any term of the
{\it graduate quotient\/} defines a {\it contravariant\/} functor from 
$\widetilde\P = \widetilde\F$ to $\O\-\mod$ admitting a {\it compatible complement\/}
[4,~Proposition~8.9]. Arguing by induction, at each step wa have to prove that some 
{\it $2\-$cocycle\/} is a {\it $2\-$coboundary\/}; the point is that we can choose a 
{\it $2\-$cocycle\/} which is {\it stable\/} with respect to all the $\widetilde\P\-$isomorphisms, and therefore we can apply the theorem above. The possibility of choosing such a 
{\it stable $2\-$cocycle\/} comes from [3,~Proposition~18.21 till~18.27].

\bigskip
\bigskip
\centerline{\large References}
\bigskip

\noindent
[1]\phantom{.} Carles Broto, Ran Levi and Bob Oliver,  {\it The homotopy theory
of fusion systems\/}, Journal of Amer. Math. Soc. 16(2003), 779-856.
 \smallskip\noindent
[2]\phantom{.} Stefan Jackowski and James McClure, {\it Homotopy
decomposition of classifying spaces via elementary abelian subgroups\/},
Topology, 31(1992), 113-132.
\smallskip\noindent
[3]\phantom{.} Llu\'\i s Puig, {\it ``Frobenius categories versus Brauer blocks''\/}, Progress in Math. 
274(2009), Birkh\"auser, Basel.
\smallskip\noindent
[4]\phantom{.} Llu\'\i s Puig, {\it Existence, uniqueness and functoriality
of the perfect locality over a Frobenius $P\-$category\/}, arxiv.org/abs/1207.0066.

\end